\documentclass{article}

\usepackage{amssymb}
\usepackage{amsthm}
\usepackage{amsmath}

\usepackage{subfig}

\usepackage{todonotes}

\newtheorem{theorem}{Theorem}[section]

\newtheorem{lemma}[theorem]{Lemma}
\newtheorem{proposition}[theorem]{Proposition}

\theoremstyle{definition}

\newtheorem{remark}[theorem]{Remark}

\newcommand{\R}{\mathbb{R}}
\newcommand{\C}{\mathbb{C}}
\newcommand{\N}{\mathbb{N}}
\renewcommand{\O}{\mathcal{O}}
\renewcommand{\P}{\mathcal{P}}
\renewcommand{\H}{H}

\newcommand{\hcont}{h_{\gamma}}	
\newcommand{\thetacont}{\theta_{\gamma}} 
\newcommand{\hcsg}{h_{\beta}}	
\newcommand{\thetacsg}{\theta_{\beta}}	

\renewcommand{\i}{\imath}
\renewcommand{\vec}{\mathbf}
\newcommand{\vecb}{\mathbf}
\newcommand{\ol}[1]{\mbox{$\overline{#1}$}}

\begin{document}

\title{GMRES-based multigrid for the complex scaled preconditoner for the indefinite Helmholtz equation}
\author{B. Reps\textsuperscript{1,2}\footnote{Corresponding author. Email: {\tt bram.reps@ua.ac.be}}, W. Vanroose\textsuperscript{1} and H. bin Zubair\textsuperscript{3}}
\date{}
\maketitle
\begin{center}
\textsuperscript{1}{\it University of Antwerp, Dept. of Mathematics and Computer Science, Middelheimlaan 1, B-2020 Antwerp, Belgium}\break
\textsuperscript{2}{\it Intel Exascience Lab, Kapeldreef 75, B-3001 Leuven, Belgium}\break
\textsuperscript{3}{\it Institute of Business Administration, Dept. of Mathematical Sciences, Faculty of Computer Science, University Rd., 75270 Karachi, Pakistan}
\end{center}

\begin{abstract}
Multigrid preconditioners and solvers for the indefinite Helmholtz equation suffer from non-stability of the stationary smoothers due to the indefinite spectrum of the operator. In this paper we explore GMRES as a replacement for the stationary smoothers of the standard multigrid method. This results in a robust and efficient solver for a complex shifted or stretched Helmholtz problem that can be used as a preconditioner. Very few GMRES iterations are required on each level to build a good multigrid method. The convergence behavior is compared to a theoretically derived stable polynomial smoother. We test this method on some benchmark problems and report on the observed convergence behavior.
\end{abstract}

\section{Introduction}
The Helmholtz equation 
\begin{equation}
\H u(\vec{x}) \equiv -(\triangle + \phi(\vec{x}))u(\vec{x}) = \chi(\vec{x})\qquad \vec{x} \in \Omega \subset \R^d \label{eqn:helm}
\end{equation}
(where $\triangle$ stands for the Laplacian) which arises as a simplified stationary model of many diverse problems, forms a fairly well-known challenge for iterative methods and has been the subject of unabated research ever since large scale electronic computations became widespread. Due to the inherent indefiniteness brought in from the underlying application, which principally occurs due to the negative wavenumber function $-\phi(\vec{x})$, traditional stationary solvers as well as their multilevel enhancements do an unsatisfactory job in efficiently finding a numerical solution for Equation~\eqref{eqn:helm}. These include all related methods prototyped by the Jacobi method as well as multigrid with standard components. The reasons of this incapacity have been analyzed in detail in many works.
Brandt and Livshits \cite{BL97,L04,LB06} introduced the wave-ray methodology where the problematic error components are factorized by a high-frequency mode and a smooth function that can be approximated on a coarser grid. In \cite{HM11} Haber and MacLachlan derive an equivalent system for \eqref{eqn:helm} by factorizing the solution itself with a Rytov decomposition. This new formulation is more feasible for numerical solvers such as e.g.\ multigrid methods. Similarly, the first-order system least-squares (FOSLS) method is based on the reformulation of second-order equations and has proven to be another succesful workaround for the arising issues with multigrid methods on indefinite Helmholtz equations \cite{CLMC94,LMCR00}.

Elman, Ernst and O'Leary \cite{EEL01} have analyzed a preconditioner for the indefinite Helmholtz equation based on a multigrid inversion of the discretized Equation~\eqref{eqn:helm} with a number of GMRES($s$) iterations used to complement the smoother. They optimized a sophisticated smoother schedule that complements Jacobi smoothing with GMRES iterations. Erlangga, Oosterlee, and Vuik \cite{EVO04,EVO06} used $\omega$-Jacobi based standard multigrid, but they applied it to the \emph{complex shifted Laplacian (CSL)} for preconditioning the original Helmholtz operator. This work is related to both of these papers. The particular indefinite linear system on which we focus, is aimed at preconditioning the outer Krylov subspace solve for \eqref{eqn:helm} through a specific multigrid inversion proposed here. This discrete preconditioner is obtained by discretizing \eqref{eqn:helm} on a complex-valued mesh. We propose complete substitution of the smoothing process with GMRES(3). This ensures that all components of the error are reduced without amplification of the smooth modes (as is the case with stationary smoothers) and results in a robust method for inverting the preconditioner with multigrid. Currently, for inverting the shifted or scaled Laplacian preconditioner, the widespread practice is the use of multigrid with under-relaxed Jacobi or ILU as a smoother \cite{EVO06, RVZ10, UMO09}. These methods may appear slightly faster in execution speed compared to the GMRES-based multigrid method proposed here. But the advantage of our choice is that all multigrid components, in contrast to ILU-based multigrid, can in principle be constructed without requiring matrix storage. Moreover, the proposed method allows a small shift size in the preconditioner that could normally only be used with ILU-based multigrid. We use and analyze this technique for the complex stretched grid preconditioner (CSG) introduced in \cite{RVZ10}.

In Section~\ref{sec:background} we provide a detailed background on the spectrum of our particular preconditioning matrix; then we test the performance of a two-grid method and a V-cycle on this matrix with the traditional smoother replaced by GMRES($s$) iterations in Section~\ref{sec:gmres_smoother}. We observe that using GMRES($3$) on all levels gives satisfactory results, however there is a wavenumber-dependent (linear) convergence rate. Our theoretical analysis of this technique is presented in Section~\ref{sec:polynomial_smoother} where a polynomial smoother of third degree is constructed. We show that the stability requirements of this polynomial smoother yields a condition on the complex parameter of our preconditioner. Finally in Section \ref{sec:polynumexp} a variety of numerical results are conducted for some 2D Helmholtz problems to benchmark the performance.

\section{Spectral prologue and iterative issues}\label{sec:background}
In this section we will briefly review some theoretical results for a simple one-dimensional Laplace model defined on a unit interval and extended with a particular absorbing layer, the \emph{exterior complex scaled (ECS)} boundary layer, mathematicaly equivalent to a Perfectly Matched Layer \cite{B94, CW94, R95}. In \cite{RVZ10} it was shown that the eigenvalues of the discrete Laplacian $L_h$ lie along a pitchfork shape figure when the Shortley-Weller finite difference scheme for non-uniform grids is used with grid distance $h\in\R$ on the interior domain $[0,1]\subset\R$ and $\hcont=he^{\i\thetacont}\in\C$ on the complex interval $[1,R_z]\subset\C$, called the ECS layer. When the discrete Helmholtz operator $H_h$ with a constant negative shift $\phi=k^2$ is considered, the spectrum is shifted westwards in the complex plane implying that the results for the Laplace model problem may be extended directly. The same extension holds for the spectrum of the CSL preconditioner $M_h^{CGL}$ where an additional vertical translation is introduced dictated by the shift size $\varepsilon$ that is employed in $(1+\varepsilon\i)k^2$.

\begin{figure}
\begin{center}
\includegraphics[width = 0.7\textwidth]{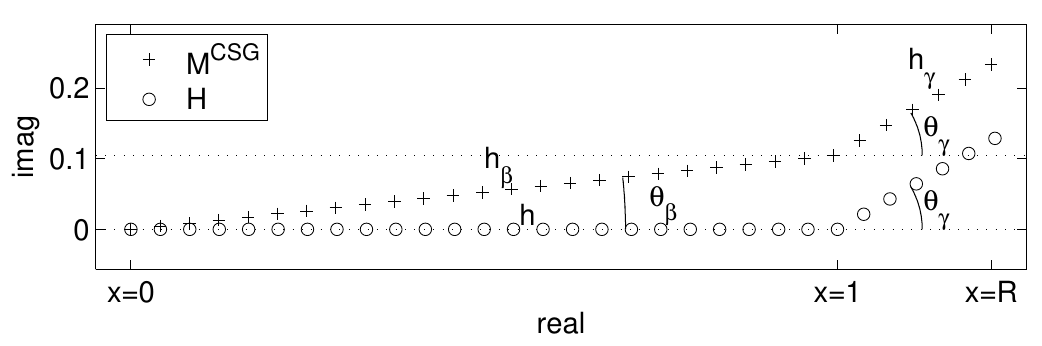}
\end{center}
\caption{One-dimensional illustration of the computational grids. The original Helmholtz operator $H$ is discretized on the ECS domain ($\circ$) and the CSG preconditioner $M^{CSG}$ on the CSG domain ($+$). The grid distance $\hcont$ on the absorbing domain extension for $x>1$ is the same for both meshes.}\label{fig:grids}
\end{figure}

Now consider the one-dimensional complex grid,
\begin{equation}\label{eq:csggrid}
(z_j)_{0\leq j\leq n+m} =
\begin{cases}
	j\hcsg, & (0\leq j\leq n);\\
	1+(j-n)\hcont,  & (n+1\leq j\leq n+m),
\end{cases}
\end{equation}
that consists of $n$ intervals of complex grid distance $\hcsg$ followed by $m$ intervals of complex grid distance $\hcont\in\C$ as illustrated in Figure~\ref{fig:grids}. It originates from an ECS grid for the discretized Helmholtz operator $H_h$ by replacing its interior real grid distance $h$ by $\hcsg = he^{\i\thetacsg}$ with $0<\thetacsg<\thetacont$. The exterior grid distance $\hcont = he^{\i\thetacont}$ provides the absorption of outgoing waves on the boundary layer. Since $\thetacont$ determines the rotation angle of the absorbing layer, we can assume that $\thetacont<\frac{\pi}{4}$ by reason of avoiding numerical reflections at the turning point between the interior and the exterior domain. Alternatively the complex grid distance can be defined as a translation in the complex plane, e.g.\ $\hcsg=h+\i h\tan\left(\thetacsg\right)$. In general we write $\hcsg = \beta h$ and $\hcont = \gamma h$ with $\beta,\gamma\in\C$.

The complex stretched grid (CSG) preconditioner $M_h^{CSG}$ is practically created by discretizing the Helmholtz operator on the grid in \eqref{eq:csggrid} with a finite difference scheme for non-uniform grids with grid distance $\hcsg\in\C$ for the interior grid points $0\leq j\leq n$ and $\hcont\in\C$ for the exterior grid points $n+1\leq j\leq n+m$. The discretization leads to a system,
\begin{equation}\label{eq:csg_dis}
\left(M_h^{CSG}\right)^{-1}H_hu_h \equiv \left[(L^{CSG}_h+k^2I_h)^{-1}(L_h+k^2I_h)\right]u_h = b_h,
\end{equation}
for the preconditioned problem. $L_h$ and $L_h^{CSG}$ represent the discretized Laplacian on the grid in \eqref{eq:csggrid} with interior rotation angle $\thetacsg=0$ and $\thetacsg>0$ respectively, $I_h$ is the identity matrix and $b_h$ is a vector containing source function values sampled at the ECS grid points. The linear systems for higher dimensional problems are created with Kronecker products of suitable one-dimensional operators and possesses similar spectral properties. Therefore, without loss of clarity, we would use the symbol $M_h^{CSG}$ interchangeably for both the one-dimensional as well as the high-dimensional analogs. Any distinction required would be stated explicitly.

\begin{lemma}\label{lem:spectriangle}
Let $M_h^{CSG}=L^{CSG}_h+k^2I_h$ be the matrix associated to the $d$-dimensional CSG preconditioner, obtained by Kronecker products of the one-dimensional Laplacian discretized on the CSG grid in \eqref{eq:csggrid}, with the Shortley-Weller finite difference scheme. Then the spectrum of $M_h^{CSG}$ is bounded in the complex plane by a triangle $\widehat{t_0 t_1 t_2}$ described by the complex points
\begin{equation*}
t_0=-k^2,\quad t_1=-k^2 + \frac{4d}{\hcsg^2} \quad \text{ and } t_2 = -k^2 + \frac{4d}{\hcont^2}.
\end{equation*}
\end{lemma}
\begin{proof}
In \cite{RVZ10} it was proven that $t_0=0,\quad t_1=\frac{4d}{\hcsg^2} \quad \text{ and } t_2 = \frac{4d}{\hcont^2}$ for the Laplacian with $d=1$ and $\thetacsg = 0$. If $d>1$, then each eigenvalue of the Laplacian $L_h^{CSG}$ is a sum of the one-dimensional operator's eigenvalues. The points $t_0, t_1$ and $t_3$ for the bounding triangle now easily follow by incorporating the negative wavenumber shift $-k^2$ in the matrix $M_h^{CSG}=L_h^{CSG}-k^2I_h$.
\end{proof}

\begin{figure}
\begin{center}
\includegraphics[width = 0.6\textwidth]{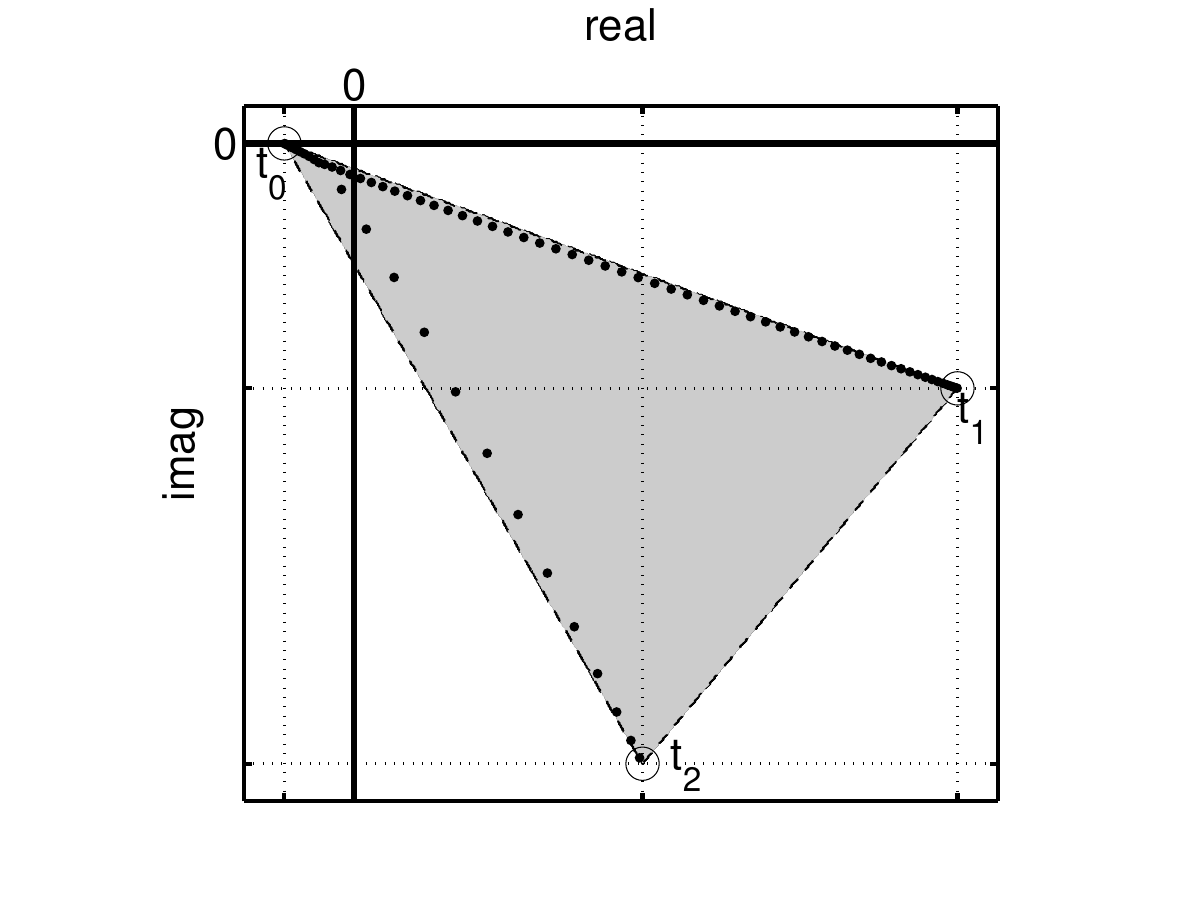}
\end{center}
\caption[Pitchfork shaped spectrum of the Helmholtz discretization matrix on an ECS grid]{The typical pitchfork shaped spectrum of the Helmholtz discretization matrix on an ECS grid. The eigenvalues lie inside a triangle region in the complex plane with vertices $t_0=-k^2$, $t_1 = -k^2+4d/\hcsg^2$ and $t_2 = -k^2+4d/\hcont^2$.}
\label{fig:pitchfork_chap_poly}
\end{figure}

\begin{figure}\label{fig:outereigvec}
\centering
\subfloat[Eigenvector $v_0$]{\includegraphics[width = 0.7\textwidth]{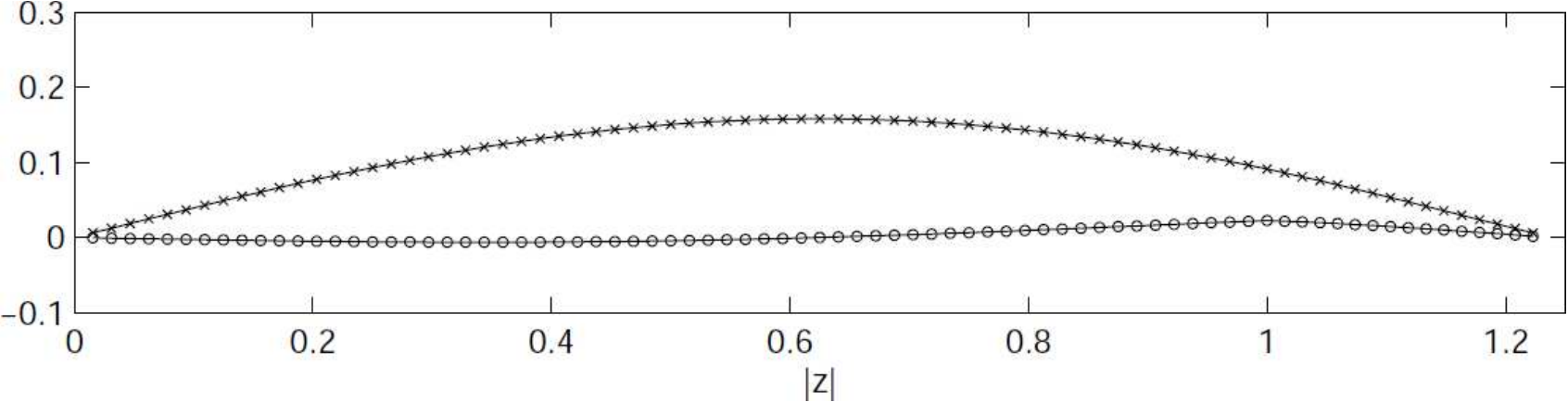}\label{fig:outereigvec_1}} \\
\subfloat[Eigenvector $v_1$]{\includegraphics[width = 0.7\textwidth]{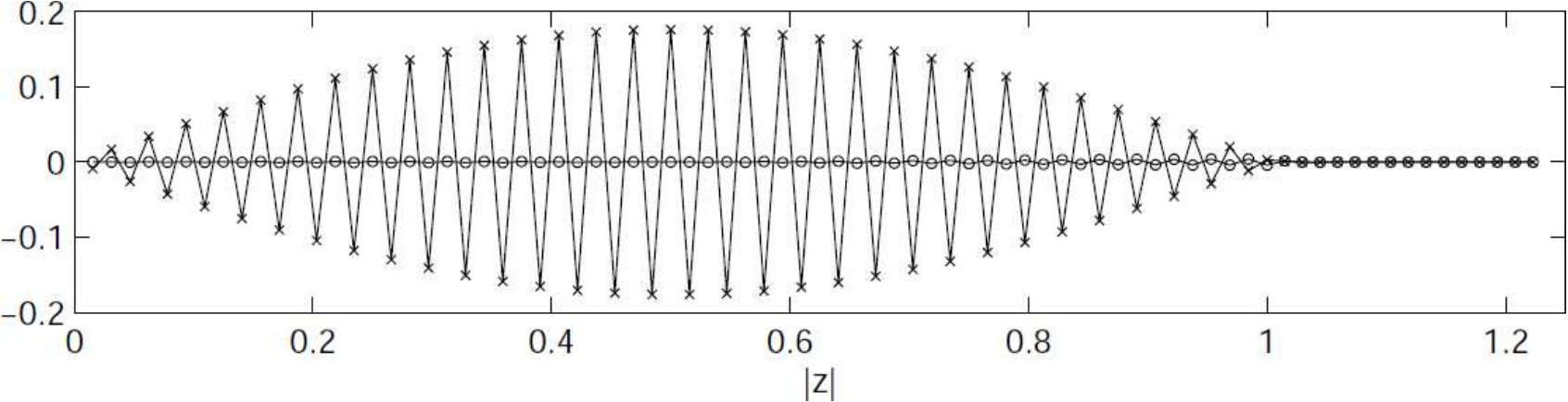}\label{fig:outereigvec_2}} \\
\subfloat[Eigenvector $v_2$]{\includegraphics[width = 0.7\textwidth]{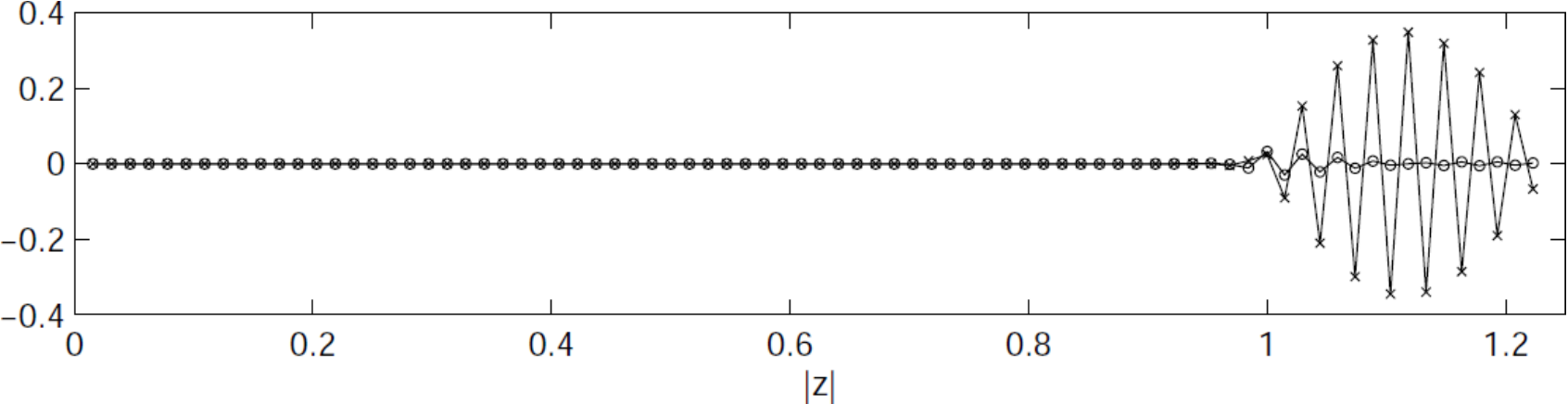}\label{fig:outereigvec_3}} 
\caption[]{Absolute value (solid line), real part (dotted line) and imaginary part (dashed line) of the three eigenvectors $v_0,v_1$ and $v_2$ associated to the extreme eigenvalues $\lambda_0\approx t_0$, $\lambda_1 \approx t_1$ and $\lambda_2 \approx t_2$ respectively. Eigenvector $v_0$ is the smoothest and is stretched over the entire domain. Eigenvector $v_1$ and $v_2$ are highly oscillatory and mainly belong to the interior and the exterior region respectively of the CSG domain.}
\end{figure}

The eigenvalues of the preconditioning matrix $M_h^{CSG}$ are not randomly distributed inside the triangle $\widehat{t_0 t_1 t_2}$, derived in Lemma~\ref{lem:spectriangle}, as their exact location is related to a physical interpretation. We discuss the one-dimensional spectrum which is illustrated in Figure~\ref{fig:pitchfork_chap_poly}.
The eigenvalues in the vicinity of the vertex $t_0=-k^2$ are aligned and correspond to the smooth or low frequency eigenvectors. We will refer to them as the \emph{smooth} eigenvalues. The eigenvector $v_0$ corresponding to the eigenvalue $\lambda_0$ closest to $t_0$ is shown in Figure~\ref{fig:outereigvec_1}. Moving further along the line the spectrum splits up into two branches at a certain point. One pronounced complex branch consists of eigenvalues with associated eigenvectors that have their largest components at indices $n\leq j\leq n+m$. Since these eigenvectors have nearly-zero components at smaller indices $1\leq j\leq n-1$, that correspond to the interior region of the grid in \eqref{eq:csggrid}, we phrase them as \emph{eigenvectors belonging to the exterior absorbing layer}. Whereas the other branch of eigenvalues in the spectrum lies closer to the real axis and corresponds to eigenvectors with their largest components at indices $1\leq j\leq n-1$, in other words, \emph{they belong to the interior domain}. Both branches in the spectrum originate in the same point on the line of smooth eigenvalues and end near the vertices $t_2 = -k^2+4d/\hcont^2$ and $t_1 = -k^2+4d/\hcsg^2$ respectively. At this end the eigenvalues correspond to highly oscillatory or high frequency eigenvectors. The eigenvector $v_1$ corresponding to the eigenvalue $\lambda_1$ that lies closest to $t_1$ is shown in Figure~\ref{fig:outereigvec_2} and eigenvector $v_2$ with $\lambda_2\approx t_2$ in Figure~\ref{fig:outereigvec_3}. Note that since the ECS grid is a special case of the CSG grid, with $\thetacsg=0$, the result in Lemma~\ref{lem:spectriangle} also holds for the original Helmholtz matrix $H_h$.

\begin{remark}
In many applications absorbing boundaries are needed on both sides of the domain. For example the grid in \eqref{eq:csggrid} is then extended on the left with an extra $m$ grid points with grid distance $\hcont$. This means that the left contour is scaled downwards in the complex plane over the same angle $\thetacont$. The spectrum of the discretized operator then has eigenvalues with algebraic multiplicity two. It is possible to choose a different ECS angle on the left contour, this will then result in an extra branch of eigenvalues in the spectrum.
When the higher dimensional Laplace operators are constructed with Kronecker products of the one-dimensional Laplacians the results on the spectrum of the discretization matrix $M^{CGS}_h$ and $H_h$ are easily extended. Every eigenvalue $\lambda$ of the $d$-dimensional Laplacian is a sum of eigenvalues $\lambda^{(i)}$ of the one-dimensional cases, $\lambda = \sum_{i=1}^{d}{\lambda^{(i)}}$. This allows for the discussion to be restricted to the corresponding one-dimensional problem without any loss of generality.
Note that actual application based Helmholtz models may require carefully engineered domains with e.g.\ smoother complex stretching or higher order discretization methods. Inasmuch as a CSG grid is used, all these generalizations can affect the internal distribution of the eigenvalues of the discretization matrix, yet the main spectral topology remains a bounded pitchfork shape with the smoothest eigenvalues aligned and the high frequency eigenvalues near $-k^2+4d/\hcont^2$ and $-k^2+4d/\hcsg^2$.
\end{remark}

\section{GMRES($s$) as a smoother substitute}\label{sec:gmres_smoother}
%
%
%
%
Multigrid methods are error correction algorithms that resolve the error on multiple grid resolutions. Their basic governing philosophy is that on a given fine resolution only particular components of the error can be reduced efficiently and that the leftover has to be addressed by switching to a coarser resolution. The main insurance of a particular multigrid method to be viable is the successful establishment of this complementarity between each pair of fine and coarse grids that may be involved in a hierarchy. In situations where this complementarity is hampered or unachievable efficiency either suffers or gives way to divergence. For this paper we assume technical familiarity with multigrid methods which are extensively described in standard literature, e.g.\ \cite{A77,BHM00,TOS01}.

Due to indefiniteness the discretized Helmholtz problem poses this complementarity issue for multigrid built with standard relaxation methods such as $\omega$-Jacobi or Gauss-Seidel. Even if the discretization matrix has all non-zero eigenvalues on the finest multigrid level one or more coarser representations can have eigenvalues undesirably close to zero that can destroy the smoothing property of the relaxation method on that level and negatively affect the entire multigrid performance. One workaround is to invest computational effort in more robust smoothers. In \cite{EEL01} classical smoothers were replaced by a sophisticated combination of GMRES and $\omega$-Jacobi on substantially indefinite levels when multigrid is applied to the Helmholtz problem, which in turn is used for preconditioning the outer GMRES solve. This requires ample GMRES iterations on certain problematic levels of the multigrid pronditioning method. In this paper we explore the idea of multigrid with the smoothing completely substituted with GMRES iterations.

\subsection{Heuristics}
We will experiment with GMRES($s$) as a substitute for the smoother on every multigrid level when it is used to invert the preconditioning matrix $M_h^{CSG}$. One smoothing step then consists of a GMRES solve on the error which is stopped after $s$ iterations, simply denoted as GMRES($s$). It constructs a minimal polynomial of order $s$ and therefore each smoothing step will cost $s$ matrix vector products. One V($\nu_1$,$\nu_2$)-cycle for example would have $\nu_1$ subsequent GMRES($s$) sweeps in the pre-smoothing phase and $\nu_2$ applications of GMRES($s$) for the post-smoothing stage. To motivate the idea of replacing the smoother by GMRES($s$) we will discuss how the spectrum of the operator evolves throughout the multigrid hierarchy.

If the operator $M_{2h}^{CSG}$ for the second finest multigrid level is built by rediscretizing $M^{CSG}=-\triangle-k^2$, then we know from Lemma~\ref{lem:spectriangle} that the spectrum is also bounded by a smaller triangle with vertices $t_0=-k^2$, $t_1=-k^2 + 4d/(2\hcsg)^2$ and $t_2 = -k^2 + 4d/(2\hcont)^2$. Again, the smooth eigenvalues are located in the neighborhood of $t_0$ and the high frequency eigenvalues near $t_1$ and $t_2$. In general, the following result holds, the proof of which is a direct consequence of Lemma~\ref{lem:spectriangle}.

\begin{proposition}\label{prop:spectriangle}
The spectrum of the rediscretized operator $M_{lh}^{CSG}$ at multigrid level $l\in\N_0$, with grid distance $l\hcsg$ in the interior domain and $l\hcont$ in the exterior, can be enclosed by the triangle $\left(\widehat{t_0 t_1 t_2}\right)_l$ with vertices
\begin{equation*}
t_0=-k^2,\quad t_1=-k^2 + \frac{4d}{l^2\hcsg^2} \quad \text{ and } t_2 = -k^2 + \frac{4d}{l^2\hcont^2}.
\end{equation*}
\end{proposition}

\begin{figure}
\subfloat[Fine level $M_h^{CSG}$]{\includegraphics[width=0.34\textwidth]{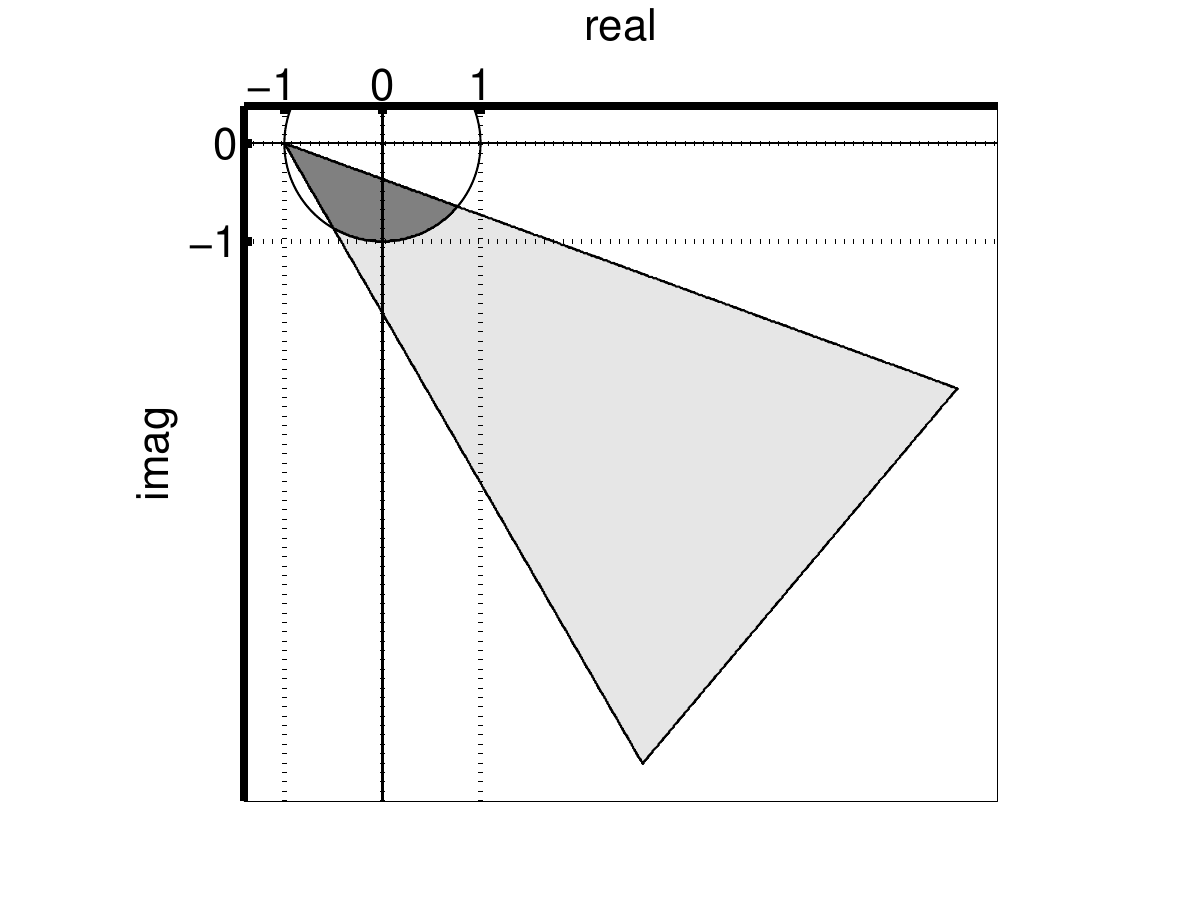}\label{fig:trianglesketch1}}
\subfloat[Intermediate level $M_{2h}^{CSG}$ ]{\includegraphics[width=0.34\textwidth]{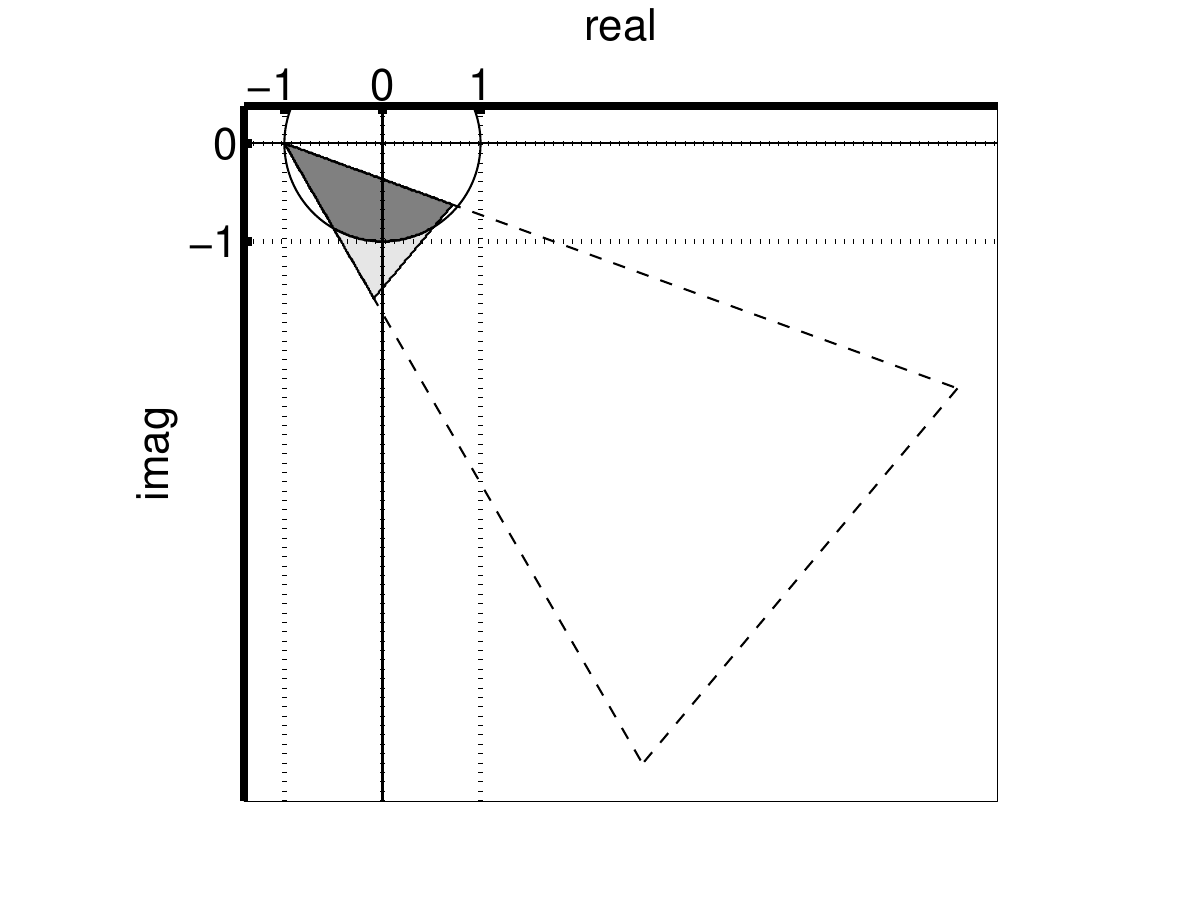}}
\subfloat[Coarse level $M_{4h}^{CSG}$ ]{\includegraphics[width=0.34\textwidth]{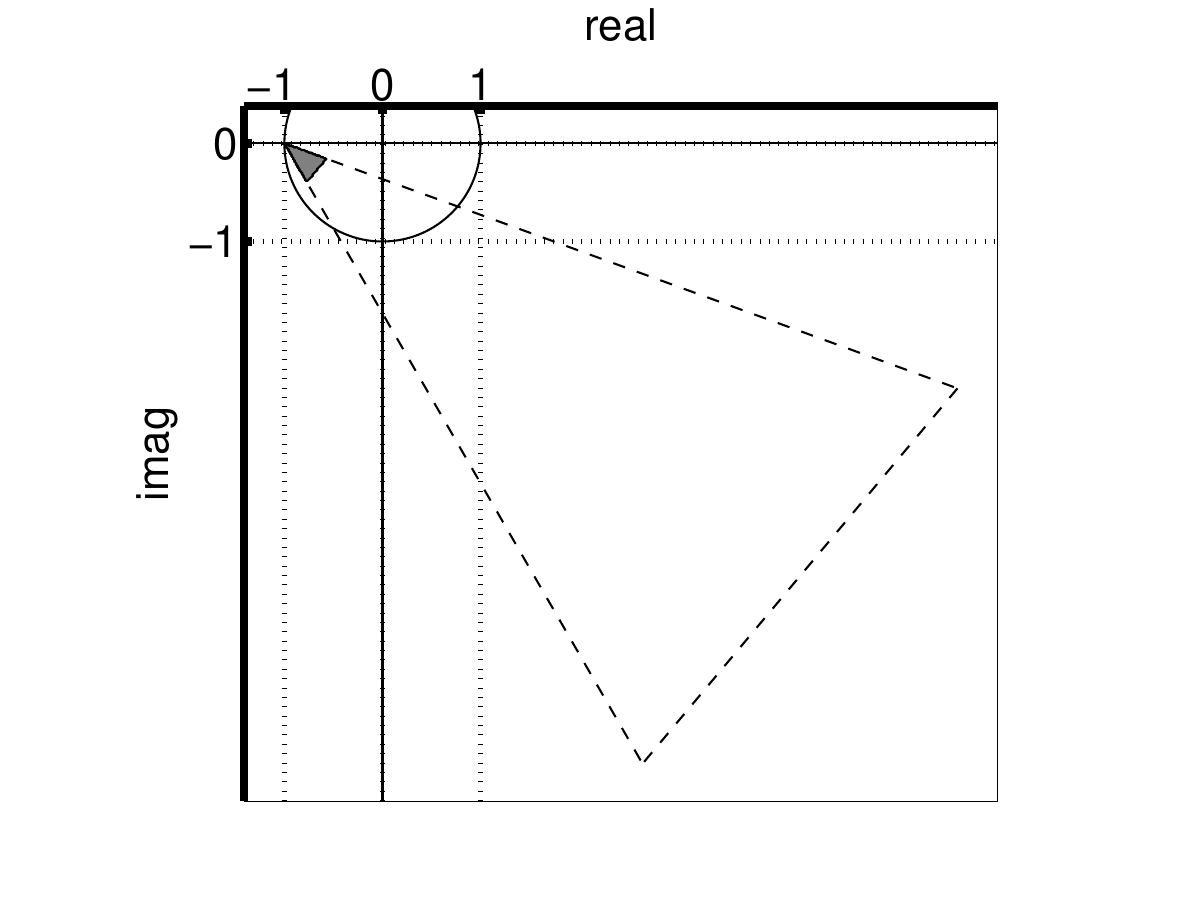}}
\caption{Bounding triangles for the spectrum of the preconditioning matrix represented on three subsequent multigrid levels by rediscretization on the coarser grids, for $k=1$ and $d=1$. The triangle of the fine level (a) with grid distances $\hcsg = he^{\i\thetacsg}$, $\hcont = he^{\i\thetacont}$ where $h=0.74$, lies mainly in the fourth quadrant of the complex plane. The spectrum of the intermediate level (b) with doubled grid distance $2h$ is equally spread over the third and fourth quadrant. On the coarsest grid (c) with grid distance $4h$ the operator is negative definite and has a small condition number.}
\label{fig:trianglesketch}
\end{figure}

In Figure~\ref{fig:trianglesketch} triangle $\left(\widehat{t_0 t_1 t_2}\right)_l$ is sketched for different levels $l$. For physically meaningful choices of the ECS rotation angle $0<\thetacont<\frac{\pi}{4}$ triangle $\left(\widehat{t_0 t_1 t_2}\right)_l$ will always lie in the lower half of the complex plane for all levels. Moreover, the upper left vertex $t_0=-k^2$ is the same for all levels and is located in the third quadrant of the complex plane. On the finest level the original real grid distance $h$ is sufficiently small in order to meet an accuracy condition such as $kh\leq0.625$ or a more stringent dependency on $k$ \cite{BGT85,IB95}. As a consequence, since in practice for the finest level $l=1$ we can assume,
\begin{equation*}
k^2\leq(0.625)^2/h^2  \ll 4d\Re\left(\frac{1}{\hcont^2}\right) \leq 4d\Re\left(\frac{1}{\hcsg^2}\right),
\end{equation*}
vertices $t_1$ and $t_2$ lie ample far into the fourth quadrant so that the largest part of triangle $\left(\widehat{t_0 t_1 t_2}\right)_l$ occupies the fourth quadrant. With vertex $t_0$ fixed at $-k^2<0$ for all multigrid levels, the width of triangle $\left(\widehat{t_0 t_1 t_2}\right)_l$ decreases with the level as $1/l^2$ and thus on the coarsest levels the triangle can be located completely in the third quadrant, i.e.\ the spectrum is negative definite.

GMRES attempts to solve a linear system $Av=b$ by minimizing the $L^2$-norm of the residual $\left\|r\right\|_2=\left\|b-Av\right\|_2$ in subsequent iterations. Let $A$ be a non-singular $N\times N$-matrix with eigenpairs $\left\{\left(\lambda_j,u_j\right) \left| \right. \left\|u_j\right\|=1 \text{ and } 1\leq j\leq N \right\}$ then the error $e$ of the current approximate solution can be written as $e = \sum_{j=1}^{N} c_j u_j$ with $c_j=\left\langle u_j,e\right\rangle$ the projection of $e$ on the space spanned by eigenvector $u_j$. For the residual this means
\begin{equation}\label{eq:residual}
  r = Ae = \sum_{j=1}^{N} \lambda_j c_j u_j,
\end{equation}
and if the eigenvectors $u_j$ are an orthonormal basis
\begin{equation*}
  \left\|r\right\|_2^2 = \left\|\sum_{j=1}^{N} \lambda_j c_j u_j\right\|_2^2 = \sum_{j=1}^{N} \left|\lambda_j\right|^2 \left|c_j\right|^2 \left\|u_j\right\|_2^2.
\end{equation*}
Suppose now that the error $e$ is not smooth and that the projections $c_j$ are of the same order for all $j=1,\ldots,N$, then the eigenvalues $\lambda_j$ with the largest magnitude will contribute the most to the residual in \eqref{eq:residual}. The GMRES method minimizes the residual and thus reduces with priority the error components in the corresponding eigenspaces of these largest eigenvalues. After the application of GMRES the projections $\hat{c}_j=\left\langle u_j,\hat{e}\right\rangle$ of the updated error $\hat{e}$ will be more biased towards the eigenvectors that have small eigenvalues.

It is clear from Proposition~\ref{prop:spectriangle} that doubling the grid distance simply reduces the range of the triangle by a factor of four, with the upper left vertex $t_0$ fixed in $-k^2$. The vertices $t_1$ and $t_2$, which are the regions of the most oscillatory eigenvalues, move closer to the smooth eigenvalues near $t_0$, while all angles in $\left(\widehat{t_0 t_1 t_2}\right)_l$ are maintained. We can now apply the above arguments on the expected performance of GMRES to the CSG preconditioner.
On the fine multigrid level in Figure~\ref{fig:trianglesketch} (a) the eigenvalues with the largest magnitude are exactly those associated to the oscillatory eigenvectors, i.e.\ points $t_1$ and $t_2$ lie far outside the circle with radius $\left|t_0\right|=k^2$, in the light grey region on the figure. After the application of GMRES the error will mainly consist of the eigenvectors that have eigenvalues in the area that lie closest to the origin, shaded in dark grey. On the fine level these are the smooth eigenvectors that can be represented on a coarser grid. This advocates the use of GMRES as a smoother substitute. On the coarse level (c) where the triangle lies completely in the third quadrant of the complex plane the condition number is significantly smaller and the problem is negative definite. Therefore, on these levels GMRES does a nice job in reducing the residual for both the smooth and oscillatory modes. Along the multigrid hierarchy we can expect a less feasible intermediate level (b). Indeed, at a certain point the triangle lies more or less equally in the third and fourth quadrant. The spectral range is still relatively large compared to case (c) and in addition there are oscillatory eigenvalues near $t_1$ which are situated in the dark grey region, meaning that their magnitude is of the same order as the smooth eigenvalues near $t_0$. However, in opposed to the original Helmholtz operator, the complex rotation of the interior region of the domain guarantees that all eigenvalues of the CSG preconditioner $M_h^{CSG}$ will be non-zero on all levels.

\subsection{Experimental observations}
The different situations described above have been confirmed through experiments with a two-grid correction scheme and a complete V-cycle. They are all run in \textsc{Matlab}\textsuperscript{\textregistered} on two quad core Intel\textsuperscript{\textregistered} Xeon CPUs (E5462 @ 2.80GHz). We focus our initial numerical experiments on solving the preconditioning problem $M_h^{CSG} u_h = b_h$ for the Helmoltz problem \eqref{eqn:helm} in 2D with constant wavenumber, $\thetacont=\frac{\pi}{6}$ for the exterior absorbing ECS regions and $\thetacsg=0.18\approx\frac{\pi}{17}$ for the interior region. This benchmark problem with right hand side $b_h$ representing a source point in the middle of the unit domain will be further evaluated in Section~\ref{sec:polynumexp}, together with more challenging Helmholtz problems. For a range of wavenumbers $0\leq k\leq180$ we measure the asymptotic convergence rates by the ratio $\|e^{(j)}\|_2/\|e^{(j-1)}\|_2$, where $e^{(j)}$ is the error in the $j$th two-grid iteration and compare GMRES($1$), GMRES($2$) and GMRES($3$) as a smoother replacement.

\begin{figure}[!ht]
\subfloat[$h=1/32$]{\includegraphics[width=0.45\textwidth]{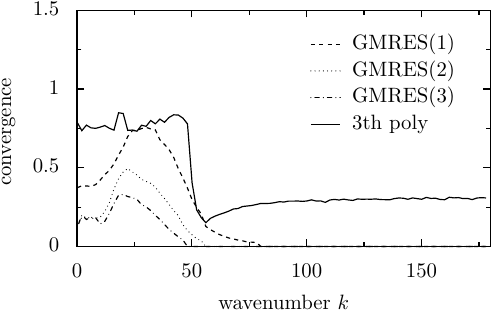}}
\subfloat[$h=1/64$]{\includegraphics[width=0.45\textwidth]{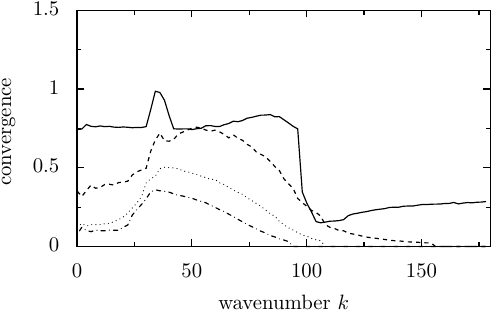}} \\
\subfloat[$h=1/128$]{\includegraphics[width=0.45\textwidth]{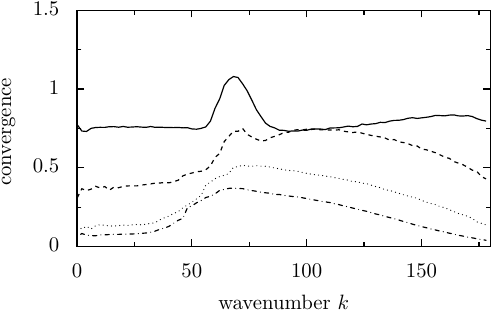}}
\subfloat[$h=1/256$]{\includegraphics[width=0.45\textwidth]{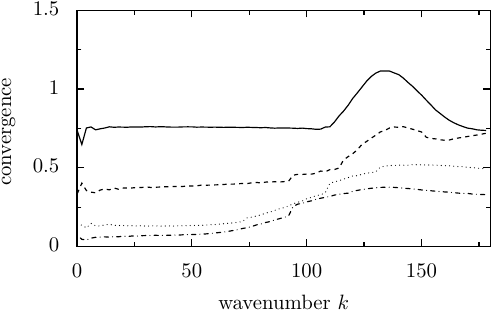}} \\
\caption{Measured two-grid convergence rate with one pre- and one post-smoothing step ($\nu_1=\nu_2=1$) for the Helmholtz equation on the CSG grid \eqref{eq:csggrid} as a function of the wavenumber $0\leq k\leq180$, for four different mesh widths $h=1/32,1/64,1/128$ and $1/256$.
  \label{fig:gmres_smoother}}
\end{figure}

This comparison is shown in Figure~\ref{fig:gmres_smoother} for the two-grid scheme for four diffent mesh widths $h=1/32,1/64,1/128$ and $1/256$. Every subfigure can be interpreted as an indicator for the action on one level of a full multigrid hierarchy. In each subfigure we can identify the same predicted behavior. First for relatively small wavenumbers the convergence rate is good because the spectrum of both the fine grid and the coarse grid operator lies mainly in the fourth quadrant and GMRES does a good job smoothing the oscillatory modes. Then, for a range of intermediate wavenumbers the spectrum of the coarse grid operator $M_{2h}^{CSG}$ lies evenly spread over the third and fourth quadrant and the two-grid operator performs clearly worse with a peak in the convergence rate. Finally for large wavenumbers both the fine and coarse grid operator are negative definite which expresses itself in a very small convergence rate in the figures for $h=1/32$ and $h=1/64$.

We conclude that for each of these two-grid tests the convergence suffers from a slowdown for a range of intermediate wavenumbers $k$ that seems to be linked to an unfavorable spectrum of the coarse grid operator. This exposes a possible weakness in a full multigrid method because for a given $k$ there will be a problematic coarser level where the operator is spread over the third and the fourth quadrant.

In the V-cycle experiments we stick to the requirent $kh\leq0.625$ to adapt the mesh width to the wavenumber. The number of cycles needed to converge to a relative residual of order $10^{-7}$ is presented in Table~\ref{tab:mgonprec} for different configurations and wavenumbers. In the experiments with the V(1,0)-cycle we see that $s=4$ does not improve the method much more than $s=3$. The experiments with $s=2$ also converge, but only after a substantially larger number of V-cycles. Increasing the number of smoothing steps ($\nu_1$,$\nu_2$) to more than (1,1) does not pay off in the eventual number of V-cycles needed to reach the tolerance. The table also shows the CPU timings for $k=80$. Based on these observations and those for the two-grid experiments, we will choose V(1,0) and V(1,1)-cycles with GMRES($3$) smoothing for the further numerical experiments presented in Section~\ref{sec:polynumexp}.

\begin{table}[!ht]
	\centering
		\begin{tabular}{ l c c c c c | c }
			\hline
			$(\nu_1,\nu_2)$ & GMRES($s$) & $k$ &    &    &    	& CPU ($k=80$)\\
											&            & 20  & 40 & 60 & 80  	&   					\\ \hline
			$(1,0)$					& $s=2$      & 18	 & 19 & 20  & 20  	& 	8.67					\\
											& $s=3$      & 13	 & 14	& 14   &15   	&	9.76 	\\
											& $s=4$      & 11	 & 12	& 12   &13  	&	12.16 	\\
							&&&&&&				\\
			$(1,1)$					& $s=2$      & 11	 & 11 & 10  & 12   & 9.14		\\
											& $s=3$      & 8	 & 9	& 9   & 9  	&	10.70		\\
											& $s=4$      & 7	 & 7	& 7   &7  		&	11.04			\\
											&&&&&&				\\
			$(2,1)$					& $s=2$      & 9	 & 9 & 8  & 10  	& 10.87 			\\
											& $s=3$      & 7	 & 7	& 7   &8  		&	13.36					\\
											& $s=4$      & 6	 & 6	& 6   & 6  	&	13.72					\\
			\hline
		\end{tabular}
	\caption{Number of multigrid V($\nu_1$,$\nu_2$)-cycles for a different number of GMRES iterations to solve the preconditioner for a constant $k$.}
	\label{tab:mgonprec}
\end{table}

\begin{figure}[!ht]
\begin{center}
\includegraphics[width=0.6\textwidth]{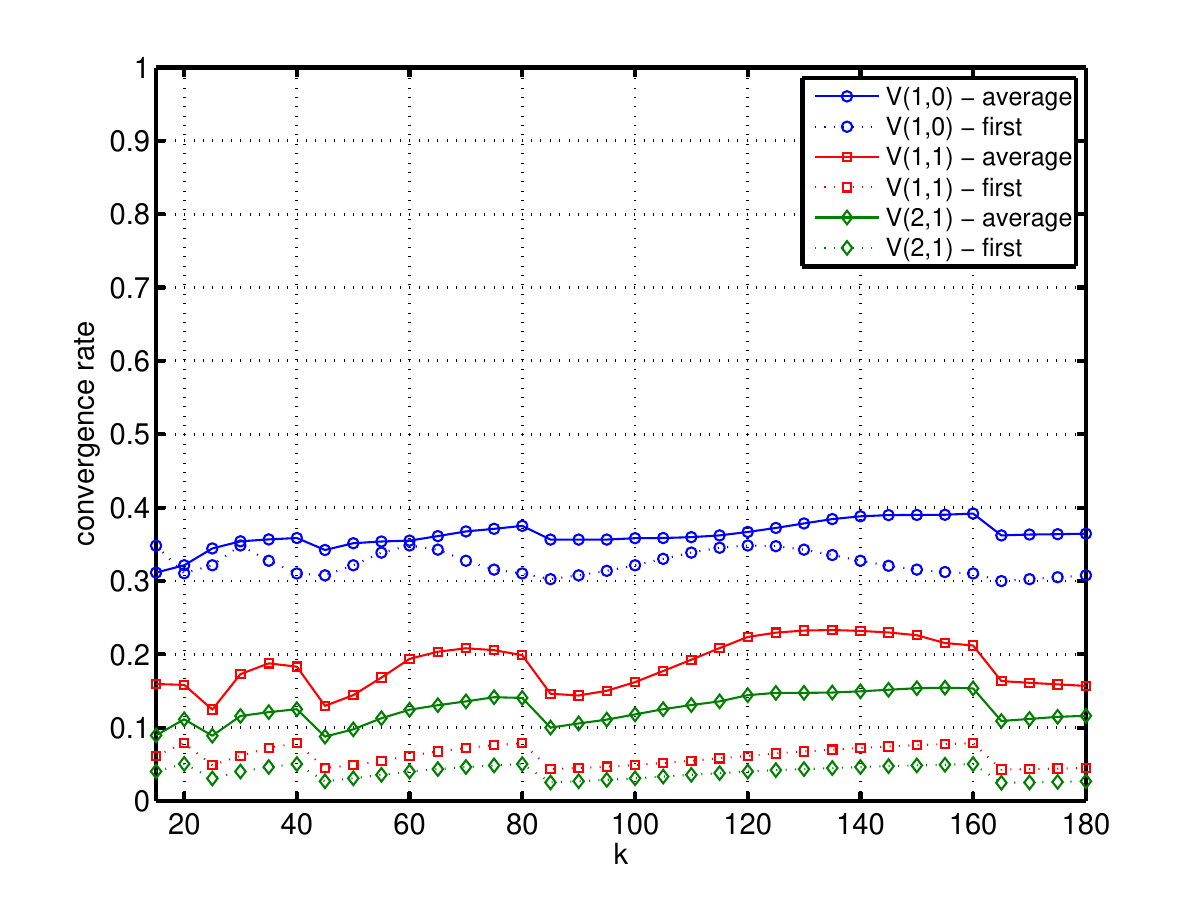}
\caption{The measured convergence rate of V(1,0) ($\circ$), V(1,1) ($\square$) and V(2,1) ($\diamond$), with GMRES($3$)-smoothing applied to the preconditioning matrix $M_h^{CSG}$ for a constant $k$ Helmoltz problem, as a function of the wavenumber $k$. The dotted lines show the residual reduction after one V-cycle, while the solid lines are average rates over more subsequent cycles.}
\label{fig:case1_k_vs_mgrate}
\end{center}
\end{figure}

In line with the two-grid experiments in Figure~\ref{fig:gmres_smoother}, the convergence rates for the different V($\nu_1$,$\nu_2$)-cycle setups with GMRES($3$) are plotted in Figure~\ref{fig:case1_k_vs_mgrate} as a function of the wavenumber $k$. The dotted lines show the residual reduction after one V-cycle, while the solid lines are average rates over more subsequent cycles. We see that the averaged convergence rate of subsequent V-cycles do not grow larger than $0.4$ for V(1,0) ($\circ$), $0.25$ for V(1,1) ($\square$) and $0.16$ for V(2,1) ($\diamond$). Note that there is a dependence on the wavenumber $k$. For just one V(1,0) and V(1,1)-cycle the convergence rate even drops under $0.35$ and $0.08$ respectively. This is an important detail because in the numerical experiments in Section~\ref{sec:polynumexp}, and by extension in practice, only one V-cycle will be used to approximately invert the preconditioning matrix $M_h^{CSG}$. The convergence rate of the V(2,1)-cycles is relatively close to that of the V(1,1)-cycles as could be expected after considering Table~\ref{tab:mgonprec}. Interestingly, these results also expose the slowdown effect that we observed in Figure~\ref{fig:gmres_smoother} for the two-grid experiments, yet for one isolated V-cycle experiment we now see several bumps in the convergence rate, instead of just one. Each individual level in the V-cycle has a specific range of wavenumbers $k$ where it does not work well and contributing to a rise in the overall multigrid convergence rate. This is especially apparent for the average rate of the V(1,1) and V(2,1)-cycles, where more relaxation steps are used per level.

\
%
\section{Explicit construction of a polynomial smoother}\label{sec:polynomial_smoother}
From the experiments in the previous section we observe that for a given grid resolution the convergence of the two-grid correction scheme suffers for a particular range of wavenumbers $k$. A similar effect was observed in the results for the V-cycle where the convergence rate showed multiple bumps. For several $k$ there is a particular problematic level of the multigrid hierarchy that causes a slowdown, seemingly related to the performance of GMRES as a smoother substitute. After three iterations the GMRES method has minimized the residual norm,
\begin{equation*}
\left\|r\right\| = \min_{p\in\P_3}\left\|p(A)r_0\right\|,
\end{equation*}
where $r_0$ is the initial residual and $A$ is the matrix that defines the linear system, by searching the space $\P_3 = \left\{p\in\C_3(x) \left|\right. p(0)=1 \right\}$ for the minimal polynomial of degree at most $3$  \cite{SS86}. In order to further understand the $k$-dependent convergence rate, we try to construct in this section a polynomial of third degree with explicit smoothing properties. We will see that there is a particular region of wavenumbers where it is harder to find such a stable polynomial smoother.

\subsection{Smoothing requirements}
Given a grid distance $h$ and wavenumber $k$, we are looking for a third order polynomial smoother $p\in \P_3$ that is intended to work for the Helmholtz problem discretized on a complex stretched grid as in \eqref{eq:csggrid}, which defines the preconditioning matrix $M_h^{CSG}$. The bounding triangle $\widehat{t_0t_1t_2}$ for the spectrum in Proposition~\ref{prop:spectriangle} will be the basis for the discussion in this section. The polynomial must have the following properties to fit the concept of smoothing.
\begin{subequations}
First, it should be stable: $\forall \lambda \in \sigma\left(M_h^{CSG}\right): \left|p\left(\lambda\right)\right| < 1$. We can ensure this by requiring that 
	\begin{equation}\label{eq:polyreq1}
	\forall t \in \widehat{t_0t_1t_2}: \left|p\left(t\right)\right| < 1.
	\end{equation}
Next, $p$ should smooth the most oscillatory eigenvectors efficiently. So we should demand that $p(t)$ maps the corresponding region of eigenvalues near $t_1 = -k^2+4d/\hcsg^2$ and $t_2 = -k^2+4d/\hcont^2$ for a $d$-dimensional problem as
	\begin{equation}\label{eq:polyreq2}
	p(t_1) = 0 \text{ and } p(t_2) = 0.
	\end{equation}
In contrast the smoother should leave the smoothest eigenvector, the one with an eigenvalue near $t_0 = -k^2$, virtually untouched. So it is required that this point is mapped to the unit circle,
\begin{equation}\label{eq:polyreq3}
 p(t_0) = p(-k^2) = e^{\i\varphi},
 \end{equation}
 with $\varphi\in\left[0,2\pi\right)$.
Finally, note that by definition $p\in\P_3$ implies
	\begin{equation}\label{eq:polyreq4}
	p(0)=1,
	\end{equation}
	which embodies the natural fixed point requirement for the exact solution of the linear system.
\end{subequations}
These conditions elaborate upon the idea of an ideal smoother: stability for all eigenmodes, the largest error reduction in the range of the most oscillatory modes and the exact solution should stay unchanged. We will show that such a polynomial smoother $p$ can be constructed if the complex stretched grid for $M_h^{CSG}$ fits certain requirements.

%
%
%
%
The polynomial can be interpreted as a sequence of three $\omega$-Jacobi steps with
different weights when written in a factored form
\begin{equation*}
p(t) = (1-\omega_1 t)(1-\omega_2 t)(1-\omega_3 t),
\end{equation*}
where $\omega_j\in\C$.

\begin{remark}
On each level of the multigrid hierarchy a different polynomial smoother $p_l\in\P_3$ must be constructed based on the bounding triangle $\left(\widehat{t_0t_1t_2}\right)_l$ in Proposition~\ref{prop:spectriangle}. Seen in this light, it is important to note that $t_0 = -k^2$ is the same for all levels and the vertices $t_1=-k^2+4d/(l\hcsg)^2$ and $t_2=-k^2+4d/(l\hcont)^2$ are level-dependent.
\end{remark}

\subsection{Construction of the polynomial}
\begin{figure}[!ht]
\begin{center}
\subfloat[Spectrum of $M_h^{CSG}$]{\includegraphics[width = 0.5\textwidth]{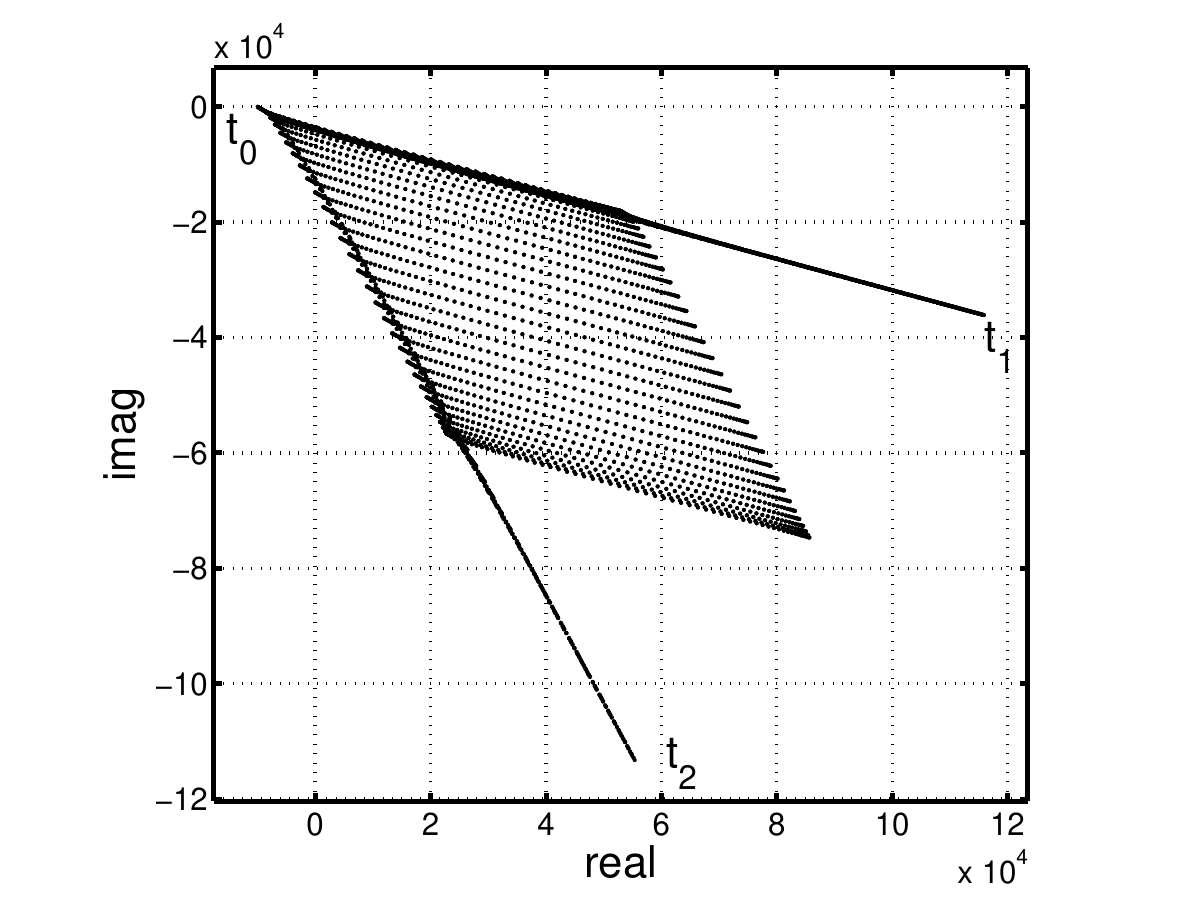}\label{sfig:spec}}
\subfloat[Polynomial map of spectrum]{\includegraphics[width = 0.5\textwidth]{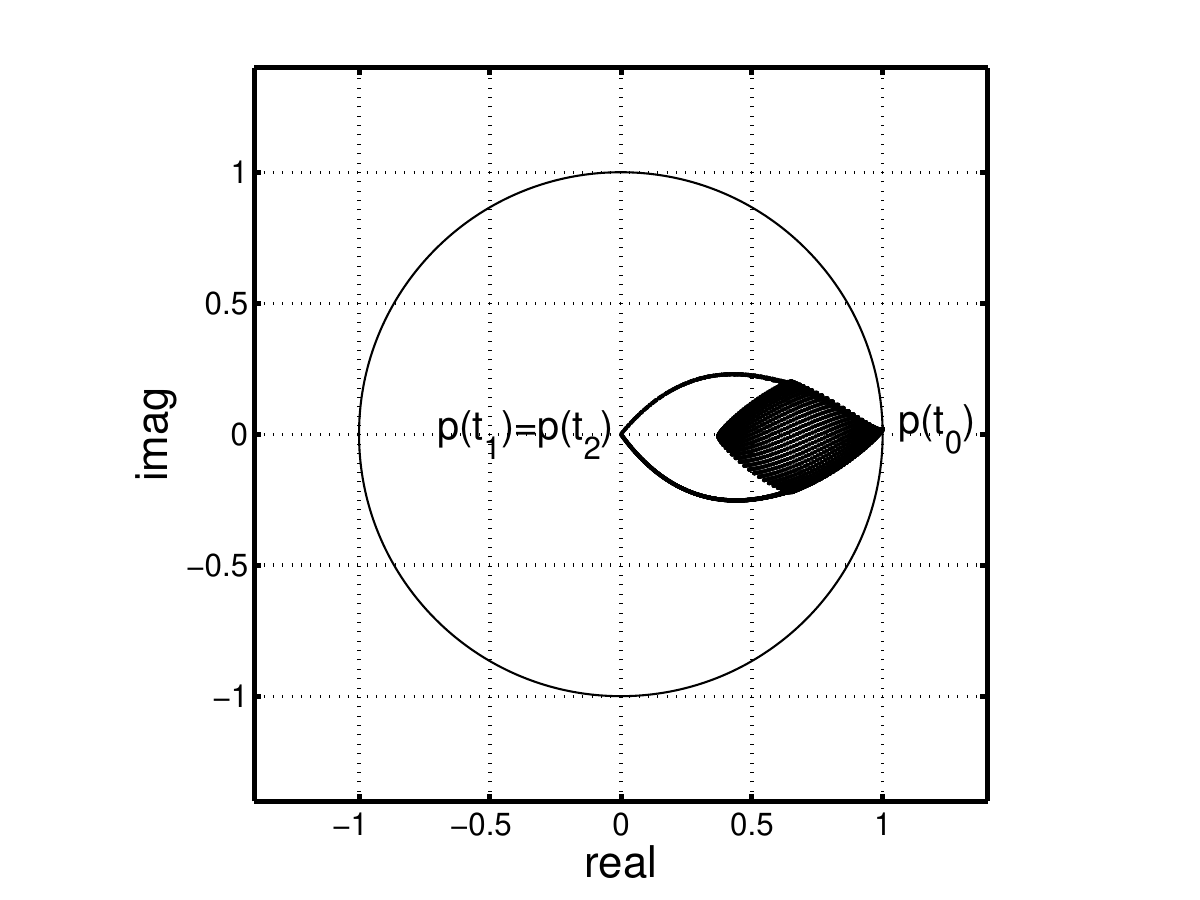}\label{sfig:full_map}}\\ 
\subfloat[Zoom of area around $p(t_0)$ for different $\varphi$]{\includegraphics[width = 0.55\textwidth]{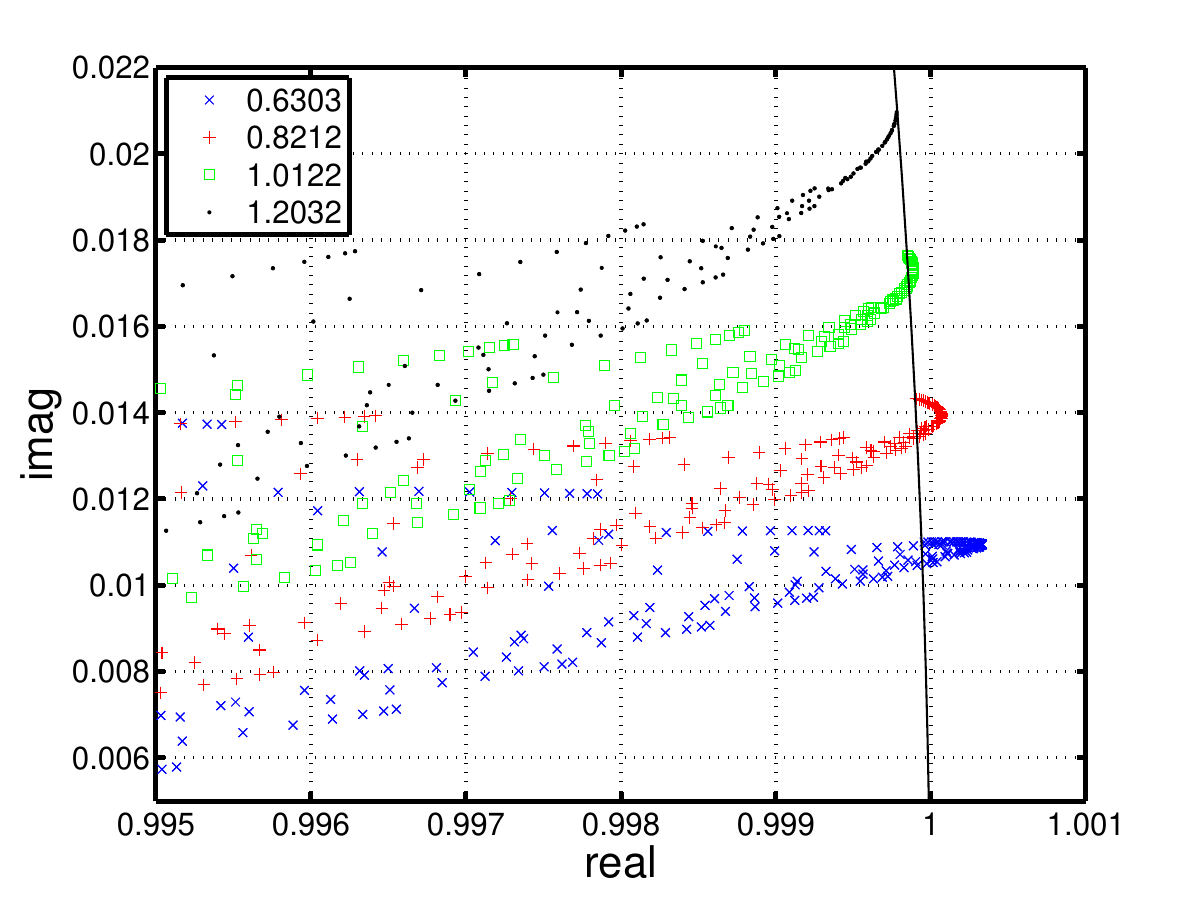}\label{sfig:detail_smooth_mode}}
\end{center}
\caption{Eigenvalues of the 2D preconditioner $M_h^{CSG}$ with constant $k=100$ on a CSG grid with $\thetacont = \pi/6$, $\thetacsg=2\pi/45$ and $h=1/128$. The spectrum is bounded by a triangle $\widehat{t_0t_1t_2}$ and can be mapped inside the unit circle by a third order polynomial $p$ with a proper choice of the parameter $\varphi$.}
\label{fig:ex1}
\end{figure}

From now on we represent $p$ in the general form,
\begin{equation}\label{eq:polynomial}
p(t) = c_0 + c_1(t-t_0) + c_2(t-t_0)^2 + c_3(t-t_0)^3,
\end{equation}
but the coefficients $c_j\in\C$ are, at this moment, unspecified. The polynomial can be viewed as a Taylor expansion around $t_0$ and so we have a simple relation between the coefficients and the derivatives
\begin{equation*}
c_0 = p(t_0),\quad c_1 = p'(t_0),\quad c_2 = p''(t_0)/2,\quad c_3 = p'''(t_0)/6.
\end{equation*}
Note that since $p\in\C_3(x)$ it is holomorphic and thus it is differentiable in the complex plane. We will first solve for the coefficients $c_0$, $c_1$, $c_2$ and $c_3$ in Equation~\eqref{eq:polynomial} such that the last three Conditions~\eqref{eq:polyreq2}-\eqref{eq:polyreq4} are met. We immediately see from Condition~\eqref{eq:polyreq3} that $c_0 = e^{\i\varphi}$. For a fixed $\varphi$ the remaining coefficients $c_1,c_2$ and $c_3$ can be calculated with polynomial interpolation. Indeed, Conditions~\eqref{eq:polyreq2} and \eqref{eq:polyreq4} then translate into the linear system,
\begin{equation} \label{eq:coefficients_system}
  \left(\begin{array}{ccc}
      -t_0 & (-t_0)^2 & (-t_0)^3 \\
      t_1-t_0 & (t_1-t_0)^2 & (t_1-t_0)^3  \\
      t_2-t_0 & (t_2-t_0)^2 & (t_2-t_0)^3  \\
    \end{array}\right)
 	\vecb{c}  = \vecb{b},
\end{equation}
with right hand side $\vecb{b} = \left(1-e^{\i\varphi}\,-e^{\i\varphi}\,-e^{\i\varphi}\right)^T$ and $\vecb{c}=\left(c_1\, c_2\, c_3\right)^T$ the vector containing the unknown coefficients. Because the rotation angle of the interior region is strictly smaller than the angle of the exterior region, $\thetacsg<\thetacont$, the three vertices of the triangle $\widehat{t_0t_1t_2}$ are distinct. It follows that the matrix in \eqref{eq:coefficients_system} is a Vandermonde matrix and so for every choice of $\varphi\in\left[0,2\pi\right)$ there is a unique solution vector,
\begin{equation}\label{eq:coefficients_solution}
\vecb{c} = \vecb{e}-e^{\i\varphi}\vecb{f},
\end{equation}
where $\vecb{e} = \left(e_1\,e_2\,e_3\right)^T$ and $\vecb{f} =\left(f_1\,f_2\,f_3\right)^T$ are the solutions of the $3\times3$ linear system \eqref{eq:coefficients_system} with right hand sides $\vecb{b} = \left(1\,\,0\,\,0\right)^T$ and $\left(1\,\,1\,\,1\right)^T$ respectively, and are independent of $\varphi$. The free parameter $\varphi$ determines the point on the unit circle where the vertex $t_0$ will be mapped by the polynomial. At this moment we have not enforced any stability constraints on $p$ yet. Not all choices of $\varphi$ lead to a stable smoother.
Indeed, the spectrum of the 2D preconditioner in Figure~\ref{sfig:spec} consists of sums of 1D eigenvalues and lies in the third and the fourth quadrant of the complex plane, bounded by the triangle $\widehat{t_0t_1t_2}$. The eigenvalues can be mapped inside the unit circle by a third order polynomial $p$ with $p(t_1) = p(t_2)=0$ and $p(t_0)=e^{\i\varphi}$ where $\varphi = 1.2032\text{ degrees}$, see Figure~\ref{sfig:full_map}. In Figure~\ref{sfig:detail_smooth_mode} we show a detailed view of the mapping around 1, for $\varphi = 0.6303$, $0.81212$, $1.0122$ and $1.2032$ degrees. Each of these angles results in a different polynomial $p$, yet not all are inside the circle.

\begin{remark}
We have made the choice to map $t_1$ and $t_2$ to zero. As a result the two most oscillatory eigenvectors will be optimally damped. This does not necessarily lead to the most efficient smoother. It might be more advantageous to map two other points inside the triangle to zero in order to obtain a better avarage damping of all the oscillatory eigenvectors. In a similar way for $\omega$-Jacobi a choice of $\omega = 2/3$ leads to a better smoother than $\omega=1/2$ for the 1D Poisson problem with Dirichlet conditions \cite{BHM00,TOS01}.
\end{remark}

\subsection{Stability condition on the parameter $\varphi$}
If we want $p$ to be a stable smoother that does not excite any of the eigenvectors then we need to map the entire triangle inside the unit circle. In particular, since the vertex in $t_0$ is mapped on the unit circle, the map of the two adjacent edges $\overline{t_0t_1}$ and $\overline{t_1t_2}$ should point inwards. Deriving an explicit sufficient condition on the mapping parameter $\varphi$ that ensures that the entire triangle is mapped inside the unit circle is technically cumbersome, yet next we will derive some useful necessary conditions that will lead to interesting insights.

We start from the requirement that the upper left corner of the triangle in $t_0$ is mapped inside the unit circle. Therefore, assume that $\varphi$ is such that $c_1= e_1+e^{\i\varphi}f_1 \neq0$, then since also $c_1=p'(t_0)$ and $p$ is holomorphic it preserves angles in $t_0$. This means we can focus on the map of only one of the two adjacent edges in $t_0$. The top edge $\overline{t_0t_1}$ of the bounding triangle, closest to the real axis, is parametrized by the line $t = t_0 + \rho e^{-\i2\thetacsg}$ with $\rho\in\left[0,\frac{4d}{h^2}\right]$ and is mapped by $p$ to
\begin{equation}\label{eq:polynomial_topedge}
p\left(t_0 + \rho e^{-\i2\thetacsg}\right) = c_0 + c_1\rho e^{-\i2\thetacsg} + c_2\rho^2 e^{-\i4\thetacsg} + c_3\rho^3 e^{-\i6\thetacsg}.
\end{equation}
We want the points on this line mapped inside the unit circle, this means
\begin{equation}\label{eq:stable_firstorder}
\left|p\left(t_0 + \rho e^{-\i2\thetacsg}\right)\right| = \left|c_0 + c_1\rho e^{-\i2\thetacsg} +R\left(\rho^2\right) \right| < 1,
\end{equation}
with $R(\rho^2)$ a correction term of order $\O(\rho^2)$. 

Because $|c_0 + c_1\rho e^{-\i2\thetacsg} +R(\rho^2)| \leq |c_0 + c_1\rho e^{-\i2\thetacsg}| +|R(\rho^2)|$, the constraint $|c_0 + c_1\rho e^{-\i2\thetacsg}|<1$ implies that the inequality in \eqref{eq:stable_firstorder} holds  for $\rho \ll 1$ sufficiently small. Using this constraint together with the solution that we derived earlier in Equation~\eqref{eq:coefficients_solution} and some basic complex algebra, allows us to derive a first condition on $\varphi$,
\begin{equation}\label{eq:firstorder_varphi}
\left|e_1\right| \cos\left(\varphi+\varphi_e-2\thetacsg\right) < |f_1|\cos\left(\varphi_f-2\thetacsg\right),
\end{equation}
where we defined $\varphi_e = \arg(e_1)$ and $\varphi_f=\arg(f_1)$ as the arguments of the first components of $\vecb{e}$ and $\vecb{f}$.
Note that $|f_1|\cos\left(\varphi_f-2\thetacsg\right)=\Re(f_1e^{-\i2\thetacsg})>0$ follows from $0<\thetacsg<\thetacont<\frac{\pi}{4}$ and as a result it is always possible to choose $\varphi$ such that the inequality in \eqref{eq:firstorder_varphi} holds by making sure that $\cos\left(\varphi+\varphi_e-2\thetacsg\right)\leq0$.
However, the range of allowed $\varphi$ is broader. If $|e_1|<|f_1|\cos\left(\varphi_f-2\thetacsg\right)$ then the inequality is always true. If instead $|e_1|\geq|f_1|\cos\left(\varphi_f-2\thetacsg\right)$, then $\varphi\in[0,2\pi)$ should be taken in the interval,
\begin{equation}\label{eq:solution_varphi}
\varphi \in \left( 2\thetacsg-\varphi_e +\varphi_0,
	2\pi +2\thetacsg-\varphi_e -\varphi_0\right) ,
\end{equation}
where $\varphi_0 = \arccos\left( \frac{\left|f_1\right|}{\left|f_1\right|}\cos\left(\varphi_f-2\thetacsg\right)\right)<\frac{\pi}{2}$ is given by the inverse cosine function that maps the interval $[-1,1]$ on $[0,\pi]$.

In other words, for all $\varphi$ in the interval \eqref{eq:solution_varphi} the corresponding polynomial $p$ maps the top adjacent edge $\ol{t_0t_1}$ of the vertex $t_0$ inwards the unit circle. Note that for the other edge $\ol{t_0t_2}$ we can derive an analogous condition, yet we make use of the property that $p$ preserves the angle between $\ol{t_0t_1}$ and $\ol{t_0t_2}$ to ensure that the image of the entire upper left corner of the bounding triangle $\widehat{t_0t_1t_2}$ points inwards the unit circle. The above condition on $\varphi$ is based on the first order approximation in Equation~\eqref{eq:stable_firstorder} of the polynomial map and is therefore not sufficient to garantuee the full stability Condition~\eqref{eq:polyreq1}.

\subsection{Second order stability condition}
\begin{figure}[!ht]
\begin{center}
\includegraphics[width=0.7\textwidth]{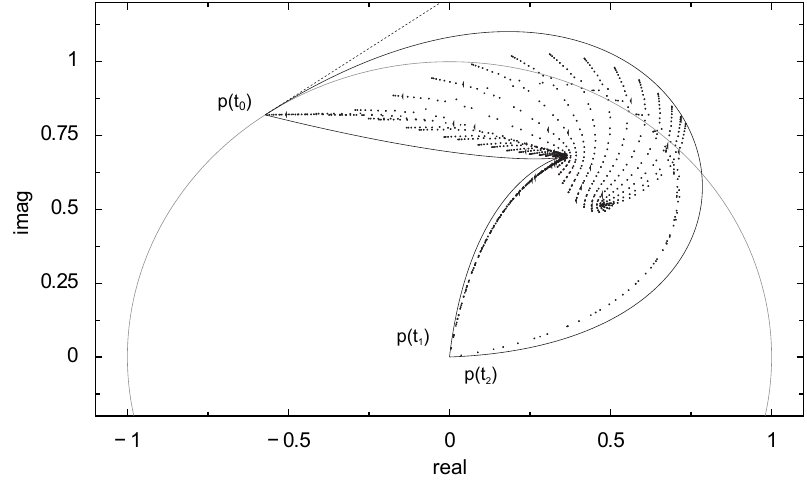}
\caption{Image of the spectrum of the preconditioner $M_h^{CSG}$ by the polynomial $p$. The first order stability condition for the parameter $\varphi$ ensures that the smooth eigenvalues are mapped inside the unit circle. Yet $p$ is unstable for some intermediate eigenvalues that lie further from the vertex $t_0$.}
\label{fig:example_unstable}
\end{center}
\end{figure}

In this section we continue the derivation of a sufficient condition on the parameter $\varphi$ for the stability requirement in \eqref{eq:polyreq1} of the polynomial smoother $p$. In Figure~\ref{fig:example_unstable} we see the action of a polynomial $p$ on the eigenvalues of the preconditioner $M_h^{CSG}$ where $\thetacsg = \frac{\pi}{36}$ and $\thetacont=\frac{\pi}{6}$, with in each dimension $n=32$ interior and $m=16$ exterior grid points.
Conditions~\eqref{eq:polyreq2}-\eqref{eq:polyreq4} are fulfilled, so $p(t_0)=e^{\i\varphi}$ lies on the unit circle and $p(t_1)=p(t_2)=0$. The first order condition for the parameter $\varphi$ in \eqref{eq:solution_varphi} is pushed to the limit. Parameter $\varphi$ is chosen such that the image of the top edge $\ol{t_0t_1}$ of the bounding triangle $\widehat{t_0t_1t_2}$ is tangent to the unit circle in the point $p(t_0)$. As a consequence, the points $t\in\ol{t_0t_1}$ near $t_0$ are mapped outside the unit circle. Nevertheless, the smooth eigenvalues are not unstable. Since $p$ preserves angles in $t_0$ and because the smooth eigenvalues of $M_h^{CSG}$ lie on a line between the two adjacent edges of this vertex, they are mapped inside the unit circle. Yet some intermediate eigenvalues at a larger distance from $t_0$ are mapped outside the circle. This illustrates how the first order stability condition for $\varphi$ is not sufficient to attain complete stability for the entire spectrum.

Next, we will derive a stricter condition on $\varphi$ for stability that is based on the second order approximation,
\begin{equation*}
\left|p\left(t_0 + \rho e^{-\i2\thetacsg}\right)\right| = \left|c_0 + c_1\rho e^{-\i2\thetacsg} + c_2\rho^2 e^{-\i4\thetacsg} + R(\rho^3)\right|<1,
\end{equation*}
of $p$ along the edge $\ol{t_0t_1}$ in Equation~\eqref{eq:polynomial_topedge}, where the correction term $R(\rho^3)$ is now of order $\O(\rho^3)$. In a similar way we impose the constraint,
\begin{equation*}
\left|c_0 + c_1\rho e^{-\i2\thetacsg} +c_2\rho^2 e^{-\i4\thetacsg}\right| < 1,
\end{equation*}
and assume that the first order stability condition $|c_0 + c_1\rho e^{-\i2\thetacsg}|<1$ is met, then this is equivalent to
\begin{equation}
2\Re\left((e^{-\i\varphi}e_2-f_2)e^{-\i4\thetacsg}\right) +\left|e_1-e^{\i\varphi}f_1\right|^2 < 0, \label{eq:secondorder_varphi}
\end{equation}
where we substituted $c_0 = e^{\i\varphi}$ and used Equation~\eqref{eq:coefficients_solution}. We arrive at a second order condition for the top edge $\ol{t_0t_1}$ to be mapped inside the unit circle. This means that it is a more stringent necessary condition on $\varphi$ than the first order condition in \eqref{eq:firstorder_varphi}. Indeed, whereas the first order condition can be fulfilled for every CSG angle $\thetacsg$ of the interior region of the domain with a proper choice of the mapping parameter $\varphi$, it is not always possible to select $\varphi$ such that the new condition in \eqref{eq:secondorder_varphi} holds. In particular if $\thetacsg$ is too small then this necessary condition cannot be met and the construction of a stable third order polynomial smoother is impossible.

\subsection{Heuristic strategy to determine the polynomial smoothers for all levels}
\begin{figure}[!ht]
\begin{center}
\includegraphics[width=0.6\textwidth]{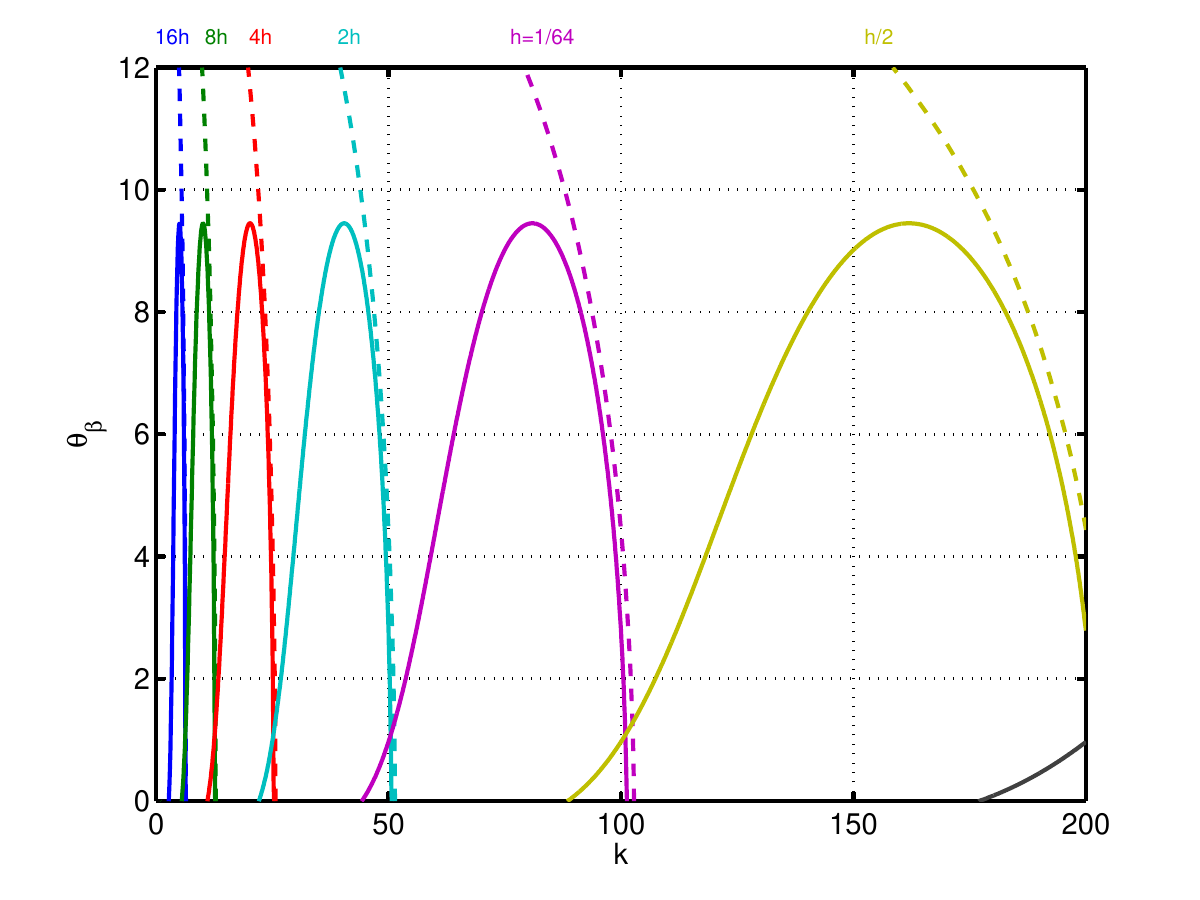}
\caption{Illustration of the first and second order necessary stability condition for $M_h^{CSG}$ with $\thetacont=\pi/6=30\text{ degrees}$, with ECS angle $\hcont=h^{\i\thetacont}$. The range of problematic wavenumbers $k$ shifts towards the left for coarser discretization. While it also narrows, the height of the region stays constant around $\thetacsg\approx9.45 \text{ degrees}$.}
\label{fig:critical_parameters}
\end{center}
\end{figure}

The two necessary conditions in \eqref{eq:firstorder_varphi} and \eqref{eq:secondorder_varphi} are illustrated in Figure~\ref{fig:critical_parameters} for a preconditioning matrix $M_h^{CSG}$ with grid distance $\hcont=h^{\i\pi/6}=30\text{ degrees}$ on the absorbing ECS layer. For a good preconditioner we want the rotation angle $\thetacsg$ of the interior grid distance $\hcsg = h^{\i\thetacsg}$ to be as small as possible. However, if $\thetacsg$ is too small then we are not able to construct a stable polynomial smoother.
The conditions are visualized in the $(k,\thetacsg)$-plane with $\thetacsg$ given in degrees for grid distances $h/2,h,2h,4h,8h,16h$ with $h=1/64$, to reflect a typical multigrid hierarchy. Each dashed line separates two regions in the plane. On the right of the line the first order condition \eqref{eq:firstorder_varphi} on $\varphi$ is always fulfilled, whereas on the left $\varphi$ must be chosen properly from the interval in \eqref{eq:solution_varphi}. The solid line further divides this left region into two. Above the solid line there is a $\varphi$ in the interval such that the second order stability condition \eqref{eq:secondorder_varphi} is met, while below the solid line it is always violated which means that there are possible unstable eigenvalues.

\begin{figure}[!ht]
\begin{center}
\subfloat[$l=6$]{\includegraphics[width = 0.33\textwidth]{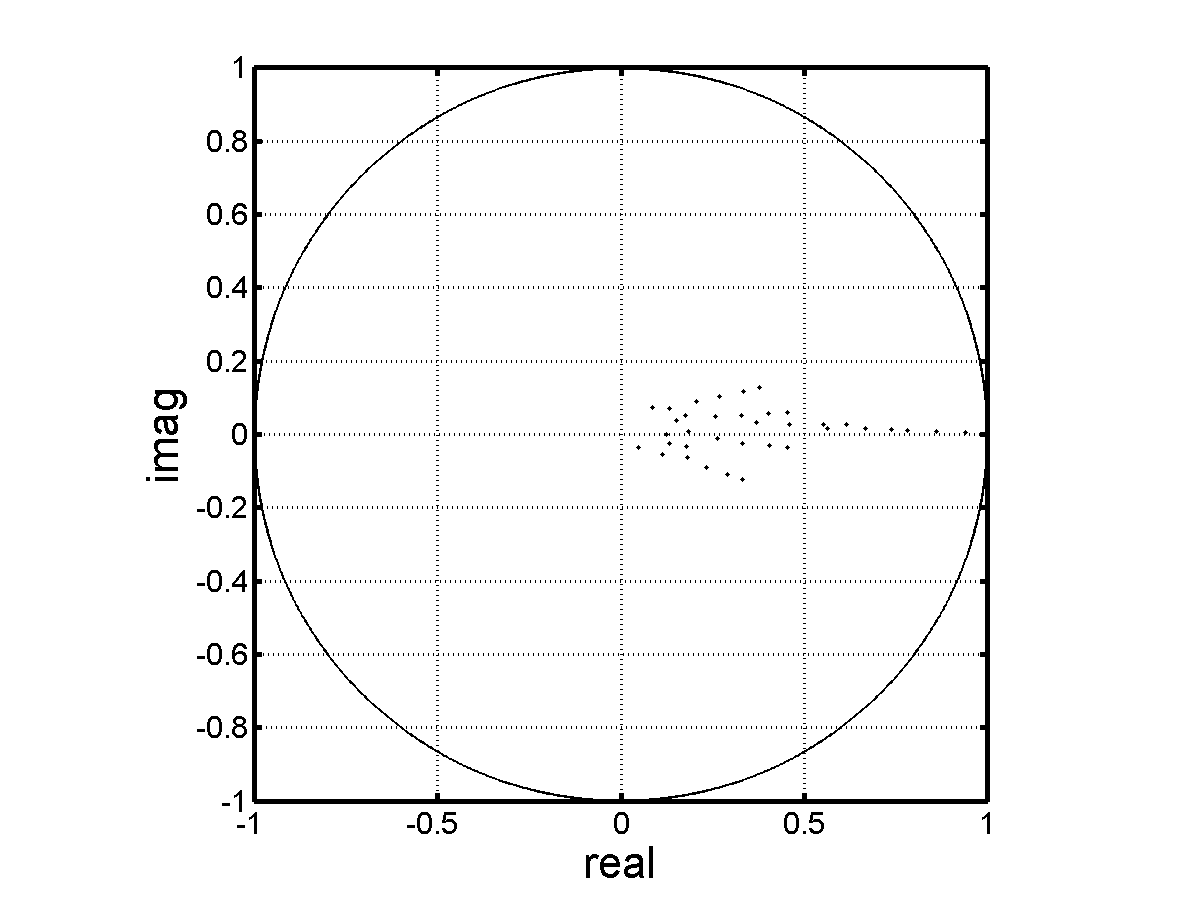}}
\subfloat[$l=5$]{\includegraphics[width = 0.33\textwidth]{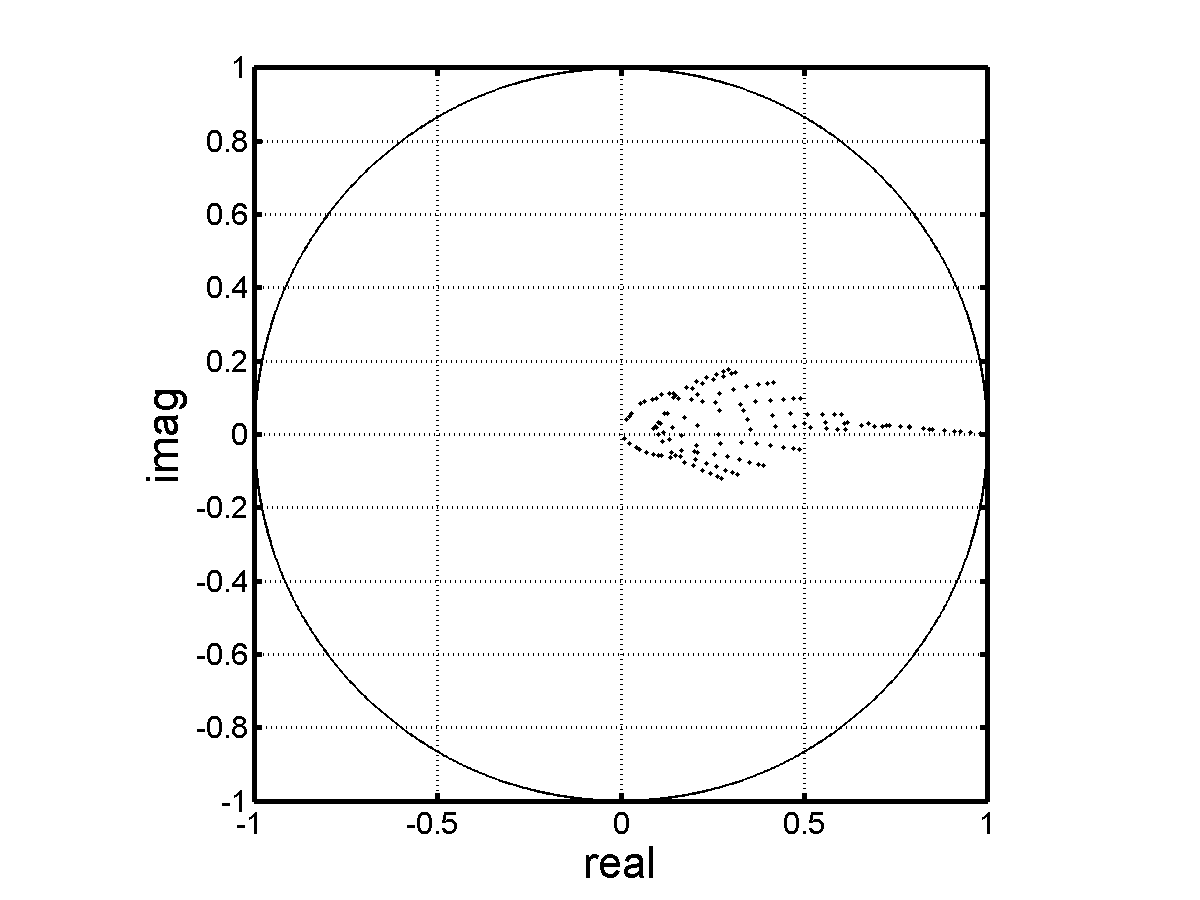}}
\subfloat[$l=4$]{\includegraphics[width = 0.33\textwidth]{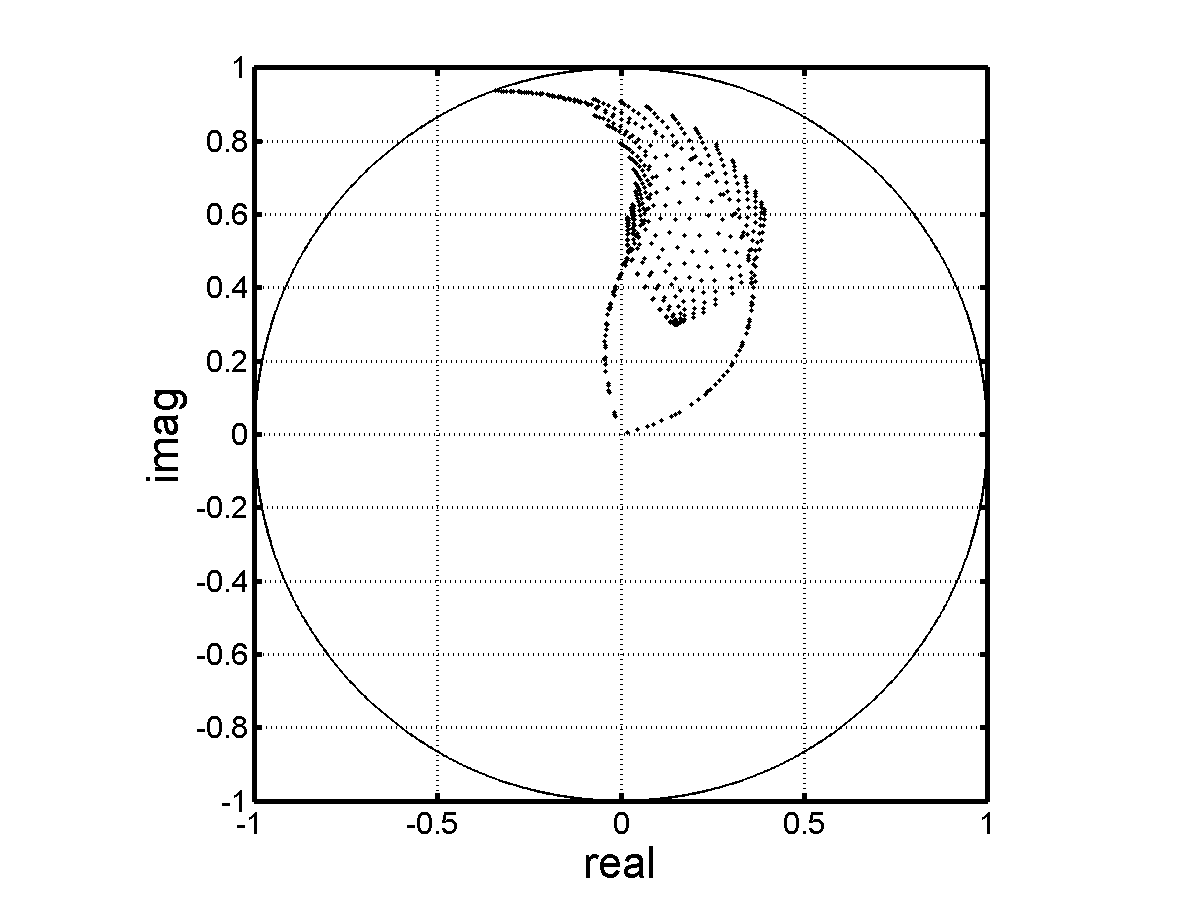}}\\
\subfloat[$l=3$]{\includegraphics[width = 0.33\textwidth]{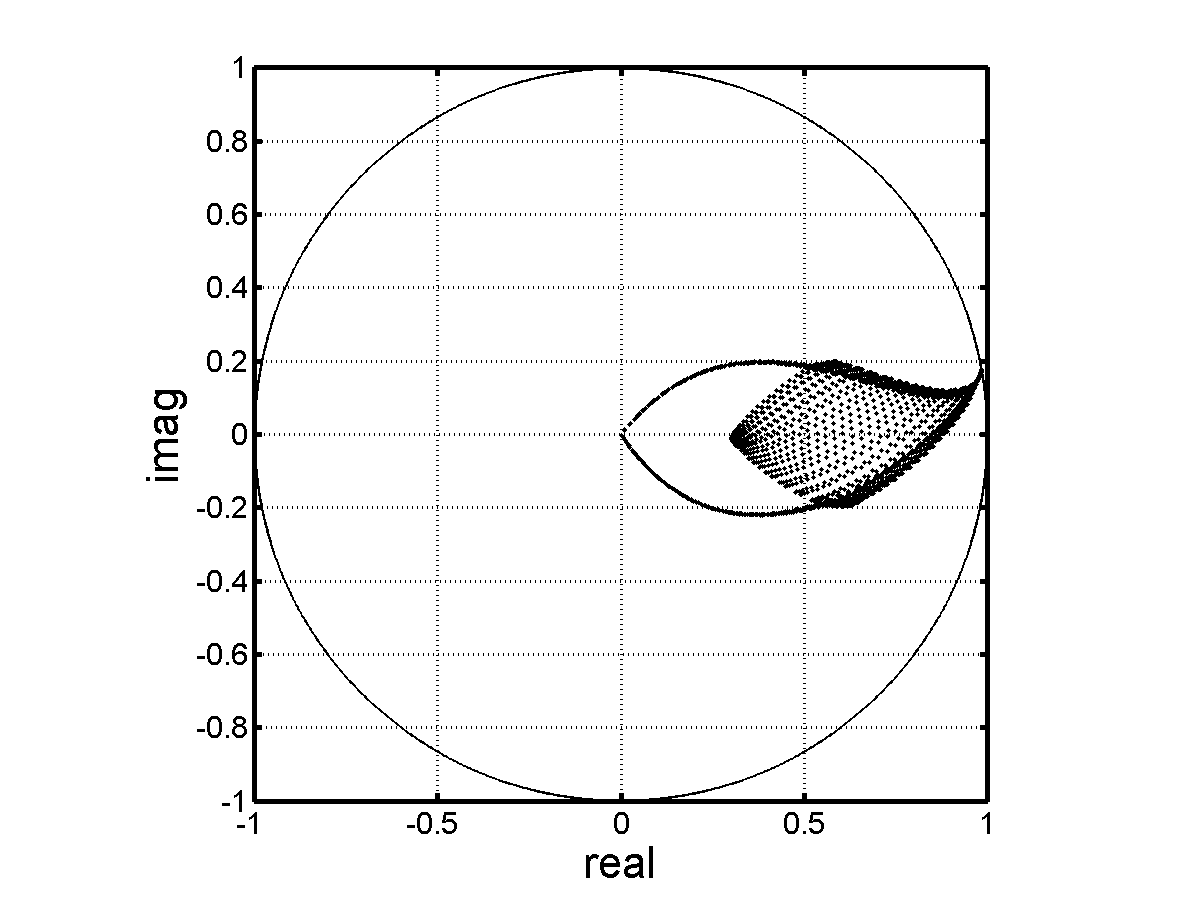}}
\subfloat[$l=2$]{\includegraphics[width = 0.33\textwidth]{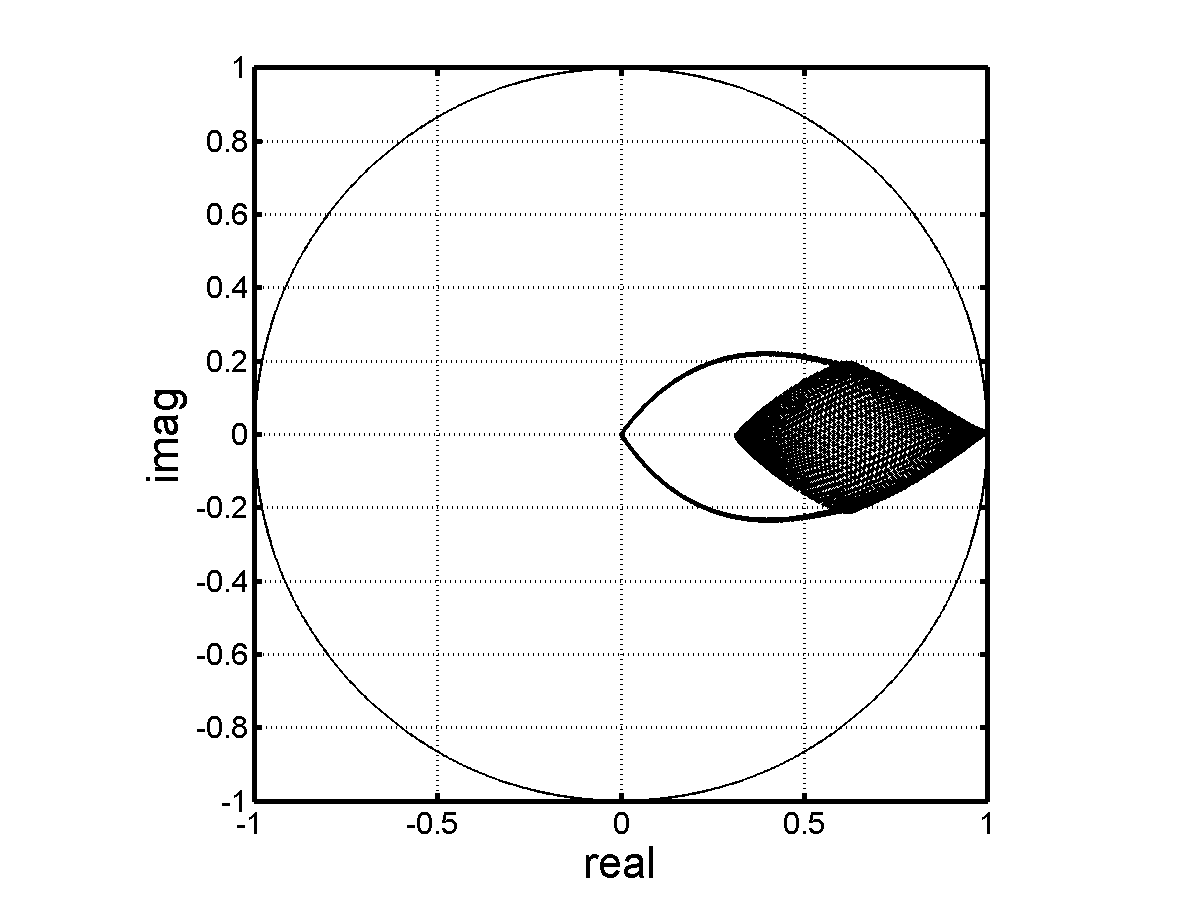}}
\subfloat[$l=1$]{\includegraphics[width = 0.33\textwidth]{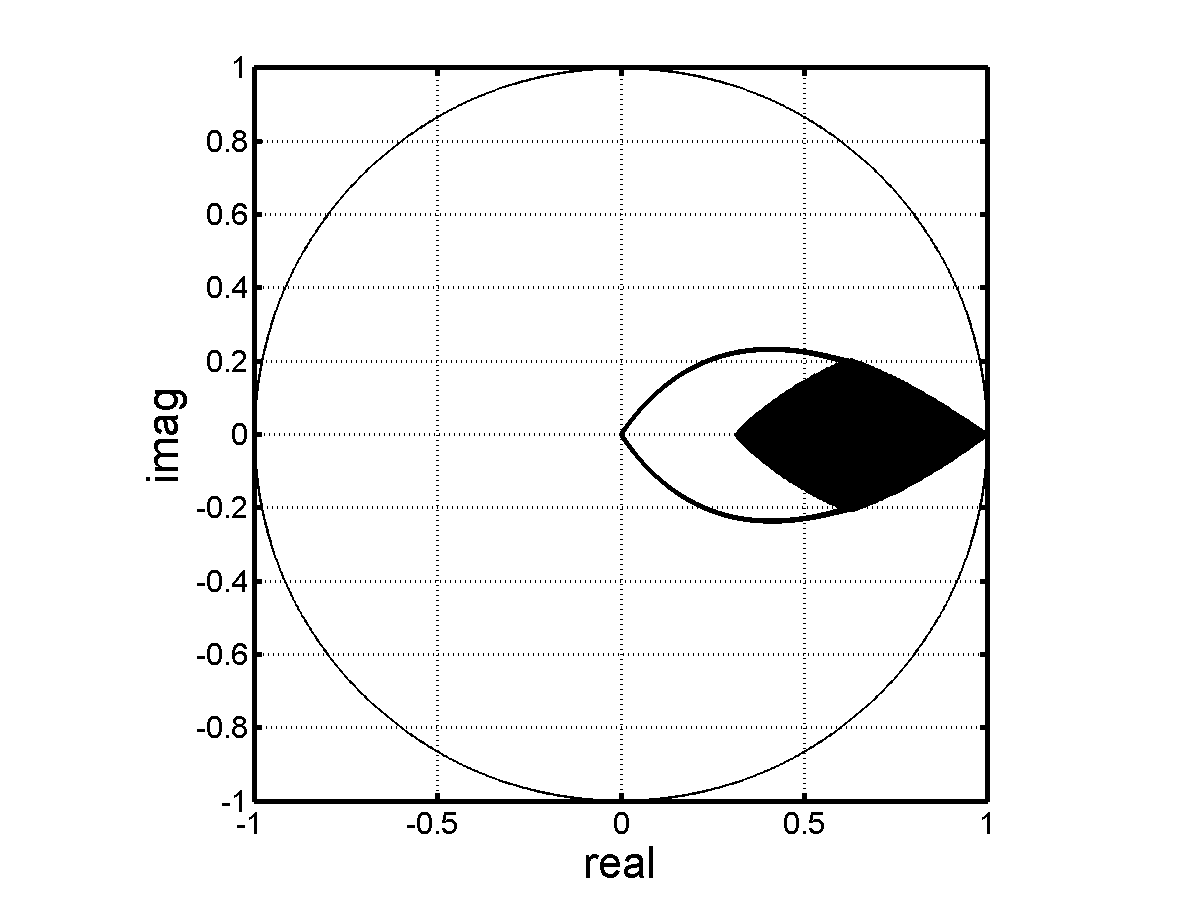}}
\caption{Image of the spectrum of preconditioning matrix $M_{lh}^{CSG}$ by the third order stable polynomials for each level $l$ of the multigrid hierarchy.}\label{fig:hierarchy_spectrum_smoother}
\end{center}
\end{figure}

In practice this means that for a given wavenumber $k$ and finest grid distance $h$ we can determine the minimal rotation angle for which there exists a stable third order polynomial smoother on every multigrid level, i.e.\ under the second order stability approximation. For the problem in this figure an angle $\thetacsg>9.5$ degrees is sufficient for all wavenumbers. For the problems we have looked at we have found that these heuristics leads to a stable smoother that maps all eigenvalues of the preconditioning matrix inside the unit circle, for every level of the multigrid hierarchy. The different maps of the spectra by the resulting polynomials are shown in Figure~\ref{fig:hierarchy_spectrum_smoother} for a 2D example with $6$ multigrid levels and $k=40$, $\thetacont=\pi/6=30\text{ degrees}$ and $\thetacsg=\pi/18=10 \text{ degrees}$. Interior mesh widths are $n=2^{8-l}$ in one dimension with additionally $m=n/4$ points on the absorbing ECS layers left and right. We see that the eigenvalues always end up inside the unit circle.

Figure~\ref{fig:critical_parameters} also illustrates that for some particular wavenumbers $k$ in the Helmholtz problem the rotation angle $\thetacsg$ of the complex stretched grid preconditioner can be very small. E.g.\ for $k\approx26$ we can choose $\thetacsg=2\text{ degrees}$. In this case $M_h^{CSG}$ and $H_h$ will be very close to each other and the outer preconditioned Krylov iteration will be highly efficient. However, for realistic problems when the wavenumber is space-dependent it is very hard to identify these wavenumbers in a robust way. It is therefore a good strategy not to choose the smallest possible $\thetacsg$, since small variations in the wavenumber could turn the smoother unstable. For this reason we choose $\thetacsg$ above the peaks of the second order stability condition and as a consequence, for each wavenumbers $k$ all the multigrid levels will have a stable smoother. However, for some wavenumbers this stability is almost pushed to the limit. With $h=1/64$ in the example, this is the case for $k\approx5,10,20,40$ and $80$. This leads to the same kind of behavior we observed in the two-grid experiments in Section~\ref{sec:gmres_smoother} with GMRES($s$) as a smoother substitute. The solid lines in Figure~\ref{fig:gmres_smoother} show the results for the polynomial smoothing method. Practically for every $k$ the convergence rate is higher than for any of the GMRES($s$) methods, with some exceptions for $s=1$. Furthermore, as a function of $k$ there is a clear sudden drop in the convergence rate for $h=1/32$ and $h=1/64$, right before $k=50$ and $k=100$ respectively. This corresponds to the dashed lines in Figure~\ref{fig:critical_parameters} that indicate unconditional stability of the polynomial smoother. Note that for the two bottom subfigures, $h=1/128$, $h=1/256$, this point does not lie in the tested wavenumber range. For wavenumbers smaller than the point determined by the dashed line, parameter $\varphi$ is selected more carefully from a level-dependent interval to ensure stability of the smoother. However, again for every grid distance $h$ there is still a region of wavenumbers where the two-grid method deteriorates, possibly due to a failing coarse grid correction, since we were able to cancel out the stability problem. However, because GMRES($s$) outperforms the polynomial smoother and since it seems to cope much better with the remaining issue of coarse grid correction, especially for $s\geq3$, we will prefer GMRES($3$) as a smoother substitute in the numerical experiments presented in the next section.


\section{Numerical experiments}\label{sec:polynumexp}
In this section we test the $CSG$ preconditioner on three benchmark problems. The experiments cover both homogeneous media with constant wavenumbers $k$ and heterogeneous media where the wavenumber is space-dependent. The discretized problem is solved with a Krylov subspace method that we will call the outer Krylov method. The complex stretched grid matrix $M_h^{CSG}$ is constructed from the same original Helmholtz equation but on a different grid with the inner part slightly rotated in the complex plane over an angle $\thetacsg=10.3132\text{ degrees}$. It is used as a preconditioner and approximately inverted with one V-cycle with GMRES($3$) as a smoother substitute. As a consequence of this non-standard smoother the actual preconditioner is not the same in every outer Krylov step. Therefore the flexible GMRES method \cite{S93} is used as outer Krylov subspace methods. The choice of the CSG angle $\thetacsg$ is based on the analysis of the previous section in the sense that the method with third order polynomial smoothing would be stable for the constant $k$ problem for all values of $k$. The number of grid points grows with $k$ according to the rule $kh\leq0.625$. We report on the number of FGMRES iterations needed to converge to a relative residual of order $10^{-7}$. The mentioned CPU times are scaled to express the computational cost per $100$ or $1000$ grid point.

%
%
\paragraph{Benchmark 1: Constant $k$ model}
In the first experiment we continue with the 2D Helmholtz problem with point source,
\begin{equation*}
\chi(x,y) = \begin{cases}
	1, \quad \text{in } x=y=\frac{1}{2};\\
	0, \quad \text{elsewhere}, \\
\end{cases}
\end{equation*}
in the middle of the unit square domain, surrounded by absorbing boundary layers. The complex shifted Laplacian preconditioner was tested on this model problem for Sommerfeld boundary conditions in \cite{EVO06} and compared to the complex stretched grid preconditioner when applied with ECS boundary conditions in \cite{RVZ10}. For the latter boundary conditions both strategies can result in a good preconditioning matrix for a Krylov subspace method but the preconditioner still needs to be (approximately) inverted by a cheap method. A multigrid method is preferred for this purpose as it is easily extended to higher dimensions where other standard methods such as ILU suffer from severe memory problems. In this paper we improve the multigrid performance on the complex stretched grid preconditioner in order to achieve a better overall convergence and to develop a robust and efficient solver for indefinite Helmholtz problems.

\begin{small}
\begin{table}[!ht]
	\centering
	\begin{footnotesize}
		\begin{tabular}{  l c c c c c }
			\hline
									$k$        		&	 20        & 40            & 60            & 80        & 100 \\
									interior grid       & $32^2$ & $64^2$   & $128^2$  & $128^2$   & $256^2$ \\ \hline							
		V(1,0)  				&  26(0.08)		& 46(0.7)  	& 54(1.65) 			& 77(2.78)  & 84(4.03) \\
	  V(1,1) 		& 17(0.08) & 26(0.48) & 34(1.27)   & 44(1.79)  & 47(2.27)  \\
		V(2,1) & 14(0.15)	&	21(2.34)		& 28(6.17)	&	37(8.34)	& 39(5.89)		 \\
			&&&&& \\
			ILU($0$)  &  66(0.18)   & 147(3.7)   & 266(20)   & 294(24.1)& 490(96.1)  \\
			ILU($0.1$)  &  77(0.20)  & 160(4.2)  & 274(22.3) & 317(30.8) & 1722(1309)  \\
			\hline
		\end{tabular}
	\end{footnotesize}
	\caption{Number of outer preconditioned GMRES iterations (and CPU time per $1000$ grid points) for different approximate inversions of $M_h^{CSG}$ for Benchmark~1. Flexible GMRES was used as outer iteration for the V-cycles with GMRES as smoother substitute.}
	\label{tab:case1_kryl}
\end{table}
\end{small}

Table~\ref{tab:case1_kryl} displays the number of iterations and the CPU time per $1000$ grid points for preconditioned FGMRES where the preconditioning matrix $M_h^{CSG}$ is approximately inverted with one V-cycle, for different values of the wavenumber $k$. The number of outer Krylov iterations with V(1,1) is significantly smaller than with V(1,0) for all values of $k$ and results in a faster CPU time. Although a V(2,1) inversion of the preconditioner requires even less outer Krylov iterations, the total CPU time is again higher than for the V(1,1) case. This confirms what has been observed in Table~\ref{tab:mgonprec} and Figure~\ref{fig:case1_k_vs_mgrate}, that extra smoothing steps do not pay off because it does not quite improve the approximate preconditioner inversion. As a reference the same problem is solved with ILU($0$) and ILU($0.1$) inversion of the preconditioner and GMRES as the outer Krylov method. Both iteration numbers and CPU times are not competitive to multigrid inversion of the preconditioner. The $k$-dependent convergence behavior is influenced by the increasing number of grid points for growing $k$ and is stronger with ILU than the V-cycle inversion. It is further studied for the CSG preconditioner in \cite{RV12}.

\begin{figure}[!ht]
\begin{center}
\includegraphics[width = 0.4\textwidth]{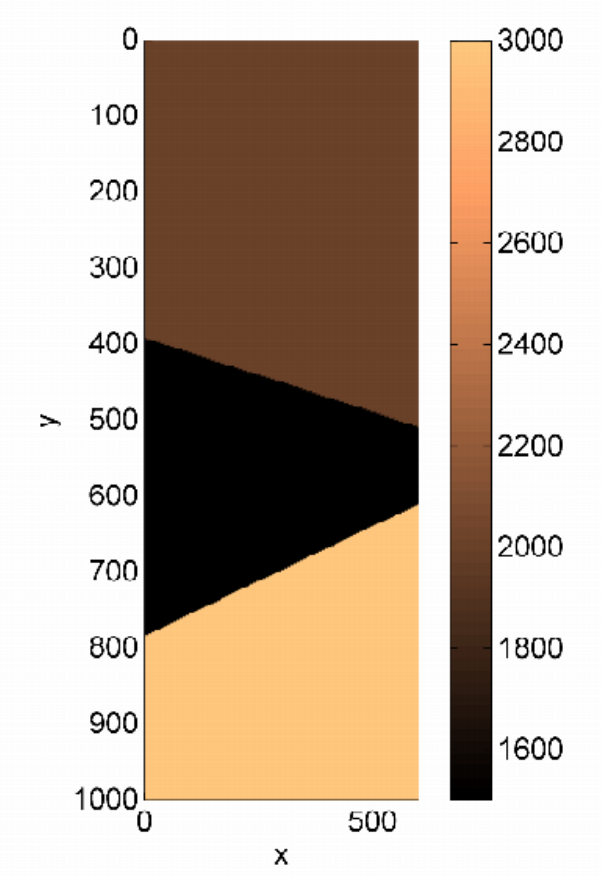}\caption{Velocity profile $c(x,y)$ for the wedge problem.}
\label{fig:wedge_velocity}
\end{center}
\end{figure}

%
%
\paragraph{Benchmark 2: Wedge model}
In the second experiment the proposed preconditioning technique is tested on a mildly heterogeneous 2D Helmholtz problem known as the wedge model, introduced in \cite{PM03} for the analysis of a preconditioner based on separation of variables and adopted in \cite{EVO06} to test the CSL preconditioner with Sommerfeld boundary conditions on all edges of the domain. In this paper the ECS boundary layers are used to absorb outgoing waves in combination with a complex stretched grid preconditioner. The rectangular domain $(0m,600m)\times(0m,1000m)$ is split in 3 regions where the speed of sound $c$ takes 3 different values (1500 m/s, 2000 m/s and 3000 m/s) as illustrated in Figure~\ref{fig:wedge_velocity}. This brings along heterogeneity in the wavenumber function that is defined as $\phi(x,y)=\left(2\pi f/c(x,y)\right)^2$, with $f\in(10Hz,50Hz)$ the frequency of the point source
\begin{equation*}
\chi(x,y) =
\begin{cases}
1,\quad \text{in } x=300,y=0;\\
0,\quad \text{elsewhere}.
\end{cases}
\end{equation*}

\begin{small}
\begin{table}[!ht]
	\centering
	\begin{footnotesize}
		\begin{tabular}{ l c c c c c }
			\hline
			          $f$      &	    10              & 20                   & 30                   & 40                    & 50 \\
			          interior grid     &$64\times128$             &$128\times256$       &$128\times256$       &$256\times512$      &$256\times512$\\
								 \hline							
	  V(1,0)  &  38(1.03)     & 68(2.57)  & 108(5.43)  & 128(9.31)& 153(12.74)  \\
		V(1,1)    &  26(0.94)     & 44(2.00)  & 68(3.56) & 79(5.15) & 96(6.95)  \\			
			\hline
		\end{tabular}
	\end{footnotesize}
	\caption{Number of outer preconditioned FGMRES iterations (and CPU time per $1000$ grid point) for V-cycle inversion of $M_h^{CSG}$ for Benchmark~2.}
	\label{tab:case2_kryl}
\end{table}
\end{small}

In Table~\ref{tab:case2_kryl} we see that V(1,1) beats V(1,0) when it comes to outer FGMRES iterations and CPU time. Other tests involved more smoothing steps and the use of ILU for approximate inversion as for Benchmark~1, but were again not competitive and are therefore left out of the table. Although the number of iterations is still rather modest, it is frequency-dependent just as the constant $k$ model.

\begin{figure}[!ht]
\begin{center}
\subfloat[Gaussian model ($k=0.4,\nu=1$).]{\includegraphics[width = 0.5\textwidth]{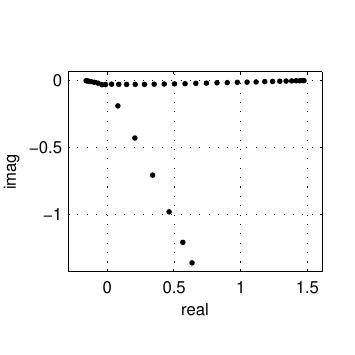}}
\subfloat[Gaussian model ($k=0.4,\nu=7$).]{\includegraphics[width = 0.5\textwidth]{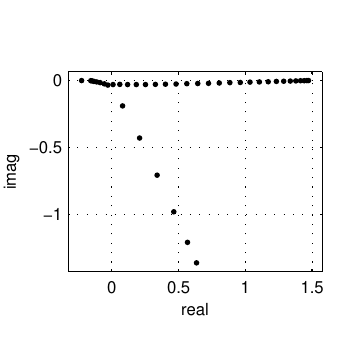}\label{fig:eig1Dgauss7}} \\
\caption{Spectra of the 1D Gaussian Helmholtz operator discretized with $32$ interior grid points.}
\label{fig:eig1D}
\end{center}
\end{figure}

%
%
\paragraph{Benchmark 3: Gaussian model}
The model for the last experiment is a 2D Helmholtz problem on the square domain $(0,50)^2$ with a strongly varying wavenumber function,
\begin{equation*}
\phi\left(x,y\right)=\nu\left(\frac{1}{e^{x^2}}+\frac{1}{e^{y^2}}\right)+k^2,
\end{equation*}
where $0<k<5$, $0<\nu<10$ and a right hand side
\begin{equation*}
\chi(x,y)=\frac{1}{e^{x^2}+e^{y^2}}.
\end{equation*}
The south and the west edges of the domain have homogeneous Dirichlet boundary conditions while ECS layers absorb outgoing waves on the north and the east edges. The problem appears in the simulation of Schr\"odinger's equation for single and multiple ionization of atoms and molecules \cite{VMRM05,CRV}. The extension to dimensions higher than 2D results in a massive amount of storage and computational complexity and is the main motivation to develop a robust matrix-free iterative method for space-dependent Helmholtz equations. It was tested for the complex shifted Laplacian and complex stretched grid preconditioners in \cite{RVZ10}. This model problem is challenging from an iterative point of view due to the highly space-dependent wave number and the parameter $\nu$. For values $\nu>2.73$ evanescent waves form near the Dirichlet edges associated to so-called bound states. These special eigenvalues of the continuous operator also appear in the spectrum of the discretized operator as one or more isolated eigenvalues on the real axis, on the left of the smoothest eigenvalues and can hamper the multigrid functionality for the preconditioning matrix, see also Figure~\ref{fig:eig1D}.

\begin{table}[!ht]
	\centering
	\begin{footnotesize}
		\begin{tabular}{  l c c c c c }
			\hline
			$k$        &	0.5  		& 1.5    		 & 2.5      & 3.5      	& 4.5 \\
interior grid       & $64^2$ 	& $128^2$ 		& $256^2$  & $512^2$ 	& $512^2$ \\ \hline							
 	V(1,0)  	& 39(0.07)  & 127(0.72) & 181(1.54) & 238(3.19) & 281(4.22) \\
	V(1,1) 	 		& 24(0.05)   & 73(0.42)   & 116(0.42) & 137(1.35)  & 194(2.40) \\
			\hline
		\end{tabular}
		\end{footnotesize}
	\caption{Number of outer preconditioned FGMRES iterations (and CPU time per $100$ grid point) for V-cycle inversion of $M_h^{CSG}$ for Benchmark~3 with $\nu=1$.}
	\label{tab:case3_kryl_mu1}
\end{table}

\begin{table}[!ht]
	\centering
	\begin{footnotesize}
		\begin{tabular}{ l c c c c c }
			\hline
		$k$        &	    0.5    & 1.5     & 2.5      & 3.5      & 4.5 \\
interior grid  & $64^2$ 	& $128^2$ 		& $256^2$  & $512^2$ 	& $512^2$ \\
								\hline
		V(1,0)  		& 139(0.45) & 159(1.67)  & 206(2.71) 	& 250(3.74) & 285(4.50) \\
		V(1,1) 	& 75(0.46)	& 106(1.52) & 134(1.67) 	& 150(1.64) & 204(2.73)  \\
			\hline
		\end{tabular}
		\end{footnotesize}
	\caption{Number of outer preconditioned FGMRES iterations (and CPU time per $100$ grid point) for V-cycle inversion of $M_h^{CSG}$ for Benchmark~3 with $\nu=7$.}
	\label{tab:case3_kryl_mu7}
\end{table}

\begin{figure}[!ht]
\begin{center}
\subfloat[$k=1.0$.]{\includegraphics[width=0.35\textwidth]{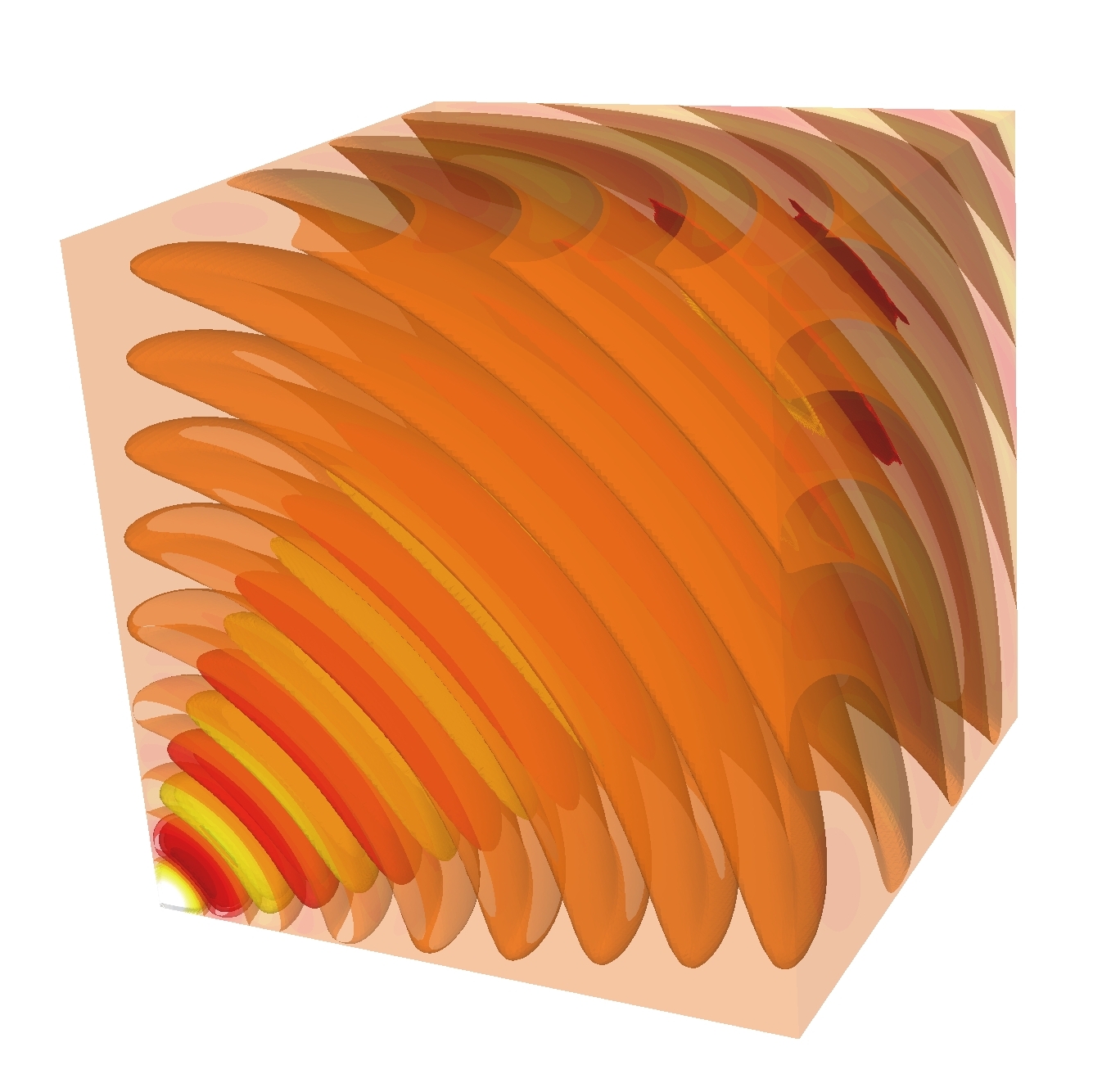}}
\hspace{0.1\textwidth}
\subfloat[$k=1.0$.]{\includegraphics[width=0.35\textwidth]{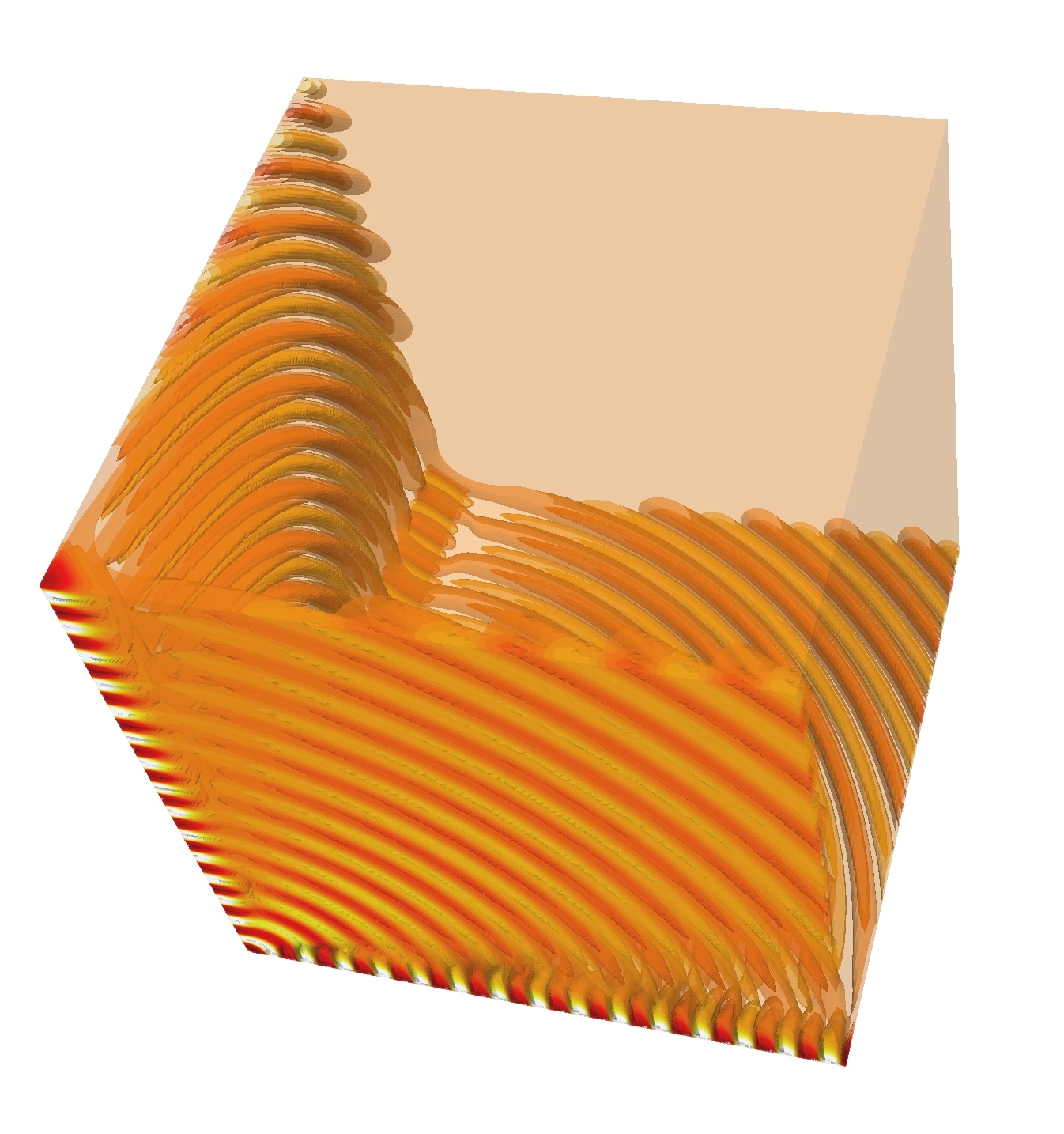}}
\caption{Real part of the solution of the 3D quantum mechanical Gaussian model with $\nu=1$ (left) and $\nu=7$ (right).}
\label{fig:gauss_3d}
\end{center}
\end{figure}

In Table~\ref{tab:case3_kryl_mu1} with the results for the model parameter $\nu=1$ we see that the method suffers again from some $k$-dependency. For a larger parameter $\nu=7$ in Table~\ref{tab:case3_kryl_mu7} we also observe this behavior. However, there is a remarkable difference in convergence speed between both problems. With $\nu=7$ the iterative solution seems harder to compute in comparison to $\nu=1$, especially for the lower wave numbers. This is most likely due to the appearance of bound states when $\nu>2.73$ as they are isolated eigenvalues in the spectrum of both the original operator and the preconditioner, positioned on the negative real axis. Especially in higher dimensions where the spectrum consists of sums of one-dimensional eigenvalues and the bound states combine with the rest of the spectrum into entire clusters of outliers they are a challenge for iterative methods.
Another effect is that they bring along the highly localized evanescent waves in the exact solution near the Dirichlet boundaries, that are clearly present in the 3D solutions for $\nu=7$ in the right panel of Figure~\ref{fig:gauss_3d}, yet absent for $\nu=1$ in the left panel.

A standard coarse grid correction in the multigrid method will not be able to solve these waves very efficiently and more complicated schemes are advisable such as L-shaped coarsening \cite{ZMO10}. It is important to note that similar results are obtained when Galerkin is used to construct the coarse operators instead of rediscretization.

\section{Conclusions and outlook}
In this paper we have analyzed the iterative solution of a Helmholtz equation discretized with finite differences and an absorbing boundary condition based on complex scaling of the domain.  The iterative solver is a flexible GMRES method that is preconditioned with a multigrid inverted complex stretched grid operator, where GRMES is used as a smoother substitute.

For each level of the multigrid hierarchy the spectrum of the discrete preconditioning operator is bounded by a triangle that lies entirely in the lower half of the complex plane. Based on the properties of the triangle we show that it is possible to choose the parameter $\thetacsg$ of the preconditioner, the rotation angle of the interior domain, such that a third order polynomial smoother is stable for all levels and all wave numbers. This smoother can be viewed as a sequence of three damped Jacobi steps with three different $\omega$'s. Note that GMRES($3$) does not necessarily have a smoothing behavior in the strict sense. However, in the experiments the convergence behavior of the GMRES based multigrid seems bounded by the results for the hand-tuned polynomial that does have the smoothing property. In addition, GMRES works in a fully automated way without tweaking. We have also observed in the various tests on model problems that the method gives satisfactory convergence results.

The numerical results are obtained for a preconditioner based on complex shifted grids but we expect that similar results will be observed for a e.g.\ complex shifted Laplacian, inasmuch the two approaches are equivalent yielding the same Krylov space convergence \cite{RVZ10}. We have found that the inversion of the preconditioner with GMRES based multigrid performs better than with pure ILU, both in outer Krylov iterations as computing time. Moreover, it requires less memory since the method is matrix-free.

In the numerical experiments we have exposed difficulties with the coarse grid correction that cannot be eliminated by the introduction of complex shifts or domain rotations only. For problems with space-dependent wave numbers that allow evanescent waves the coarse grid correction can still be problematic. Because these strongly localized waves are bound to isolated eigenvalues the use of a deflated Krylov method might be more effective which is a possible subject for future research. Another interesting outlook is the efficient implementation for 3D problems. On modern hardware stencil computations are typically communication bound since there are only a few floating point operations for each read from the slow memory. For a polynomial smoother however, communication avoiding optimizations are possible that increase the number of floating point operations per load, which have lesser overhead than communication.

\section*{Acknowledgement}
This research was partially supported by FWO-Flanders through grant G.0174.08, by a starting grant from the University of Antwerp, Belgium, by Intel\textsuperscript{\textregistered} and by the Institute for the Promotion of Innovation through Science and Technology in Flanders (IWT).\\

\bibliographystyle{unsrt}
\bibliography{main}
\end{document}